\documentclass{amsart}
\usepackage{amsmath,amssymb}
  \usepackage{paralist}
  \usepackage{graphics} 
  \usepackage{epsfig} 
\usepackage{graphicx}  \usepackage{epstopdf}
 \usepackage[colorlinks=true]{hyperref}
\hypersetup{urlcolor=blue, citecolor=red}
\usepackage{hyperref}

\usepackage{rotating}
\usepackage{adjustbox}
\usepackage{array}      
\usepackage{multicol}        
\usepackage[bottom]{footmisc} 
\usepackage{float}
\usepackage{cancel}
\usepackage[usenames,dvipsnames]{pstricks}
\usepackage{pst-grad} 
\usepackage{pst-plot} 
\usepackage{multirow} 
\usepackage{color}
\usepackage{setspace}
\usepackage{colortbl}
\usepackage{caption}

\newcommand{\C}{\mathcal{C}}
\newcommand{\R}{\mathbb{R}}

\newcommand{\wh}{\widehat}
\newcommand{\f}[2]{\frac{\displaystyle{#1}}{\displaystyle{#2}}}
\newcommand{\n}{\newline}

\def\be{\begin{equation}}
\def\ee{\end{equation}}

\def\R{\mathbb{R}}

\def\C{\mathcal{C}}

\def\f{\frac}

\def\ez{e^{\alpha z^*}}
\def\ezm{e^{-\alpha z^*}}
\def\ezo{e^{\alpha z^*(\theta_t\omega)}}

\def\azbn{\alpha z^*_{\beta,\delta}(\theta_t\omega)}
\def\azbns{\alpha z^*_{\beta,\delta}(\theta_s\omega)}
\def\ben{\begin{eqnarray}}
\def\een{\end{eqnarray}}
\def\ee{e^{\alpha z^*}}

\allowdisplaybreaks

\newtheorem{theorem}{Theorem}[section]

\newtheorem{proposition}{Proposition}[section]

\theoremstyle{definition}

\newtheorem{remark}{Remark}

\title[Random chemostat models with wall growth] 
      {Random and stochastic disturbances on the input flow in chemostat models with wall growth}

\author[J. L\'opez-de-la-Cruz]{}


 \email{jlopez78@us.es}

\thanks{Partially supported by Junta de Andaluc\'ia under Proyecto de Excelencia P12-FQM-1492.}

\thanks{$^*$ Corresponding author: jlopez78@us.es}

\begin{document}


\maketitle



\centerline{\scshape Javier L\'opez-de-la-Cruz$^*$}
\medskip
{\footnotesize
 \centerline{Dpto. Ecuaciones Diferenciales y An\'alisis Num\'erico}
   \centerline{Facultad de Matem\'aticas, Universidad de Sevilla}
   \centerline{C/ Tarfia s/n, Sevilla 41012, Spain}
} 

%

\bigskip


\begin{abstract}
In this paper we analyze a chemostat model with wall growth where the input flow is affected by two different stochastic processes: the well-known standard Wiener process, which leads into several drawbacks from the biological point of view, and a suitable Orsntein-Uhlenbeck process depending on some parameters which allow us to control the noise to be bounded inside some interval that can be fixed previously by practitioners. Thanks to this last approach, which has already proved to be very realistic when modeling other simplest chemostat models, it will be possible to prove the persistence and coexistence of the species in the model without needing the theory of random dynamical systems and pullback attractors needed when dealing with the Wiener process which, moreover, does not provide much information about the long-time behavior of the systems in many situations.
\end{abstract}

\section{Introduction}\label{intro}

The chemostat is a bioreactor very used to study both biological and ecological processes and has many applications in the real life, for instance, it is very useful when modeling and analyzing wastewater treatment processes, genetically altered organisms, antibiotic production models and fermentation models of beer and wine, to name some of the most interesting applications.\n

It is worth mentioning that the chemostat has been subjected to a large number of scientific publications and books, not only in biology and ecology but also in mathematics, in fact, there exists a recent area which is call {\it the theory of chemostat} where many researchers are involved in the last years. This interest arises due to the fact that the chemostat device can be modeled mathematically in a very simple way which reproduces very well the real process and this is the reason which encourages us to study this model.\n

\begin{figure}[H]
\begin{center}
\includegraphics[width=0.47\textwidth]{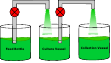}
\end{center}
\caption{The chemostat model}
\label{s}
\end{figure}

The simplest chemostat consists of three tanks, the {\it feed bottle}, the {\it culture vessel} and the {\it collection vessel}, which are connected by pumps as in Figure \ref{s}. The substrate or nutrient is stored in the feed bottle and supplied to the culture vessel, where the interactions between the substrate and the microorganisms (which are also called species) take place. Apart from this, in order to keep the volumen in the culture vessel constant, another flow is pumped from there to the collection vessel. In addition, it is also known that microorganisms tend to adhere to the walls of the culture vessel hence we would have to consider two different species in the culture vessel, the ones which are in the liquid media and those ones sticked on to the walls of the recipient.\n

Then, the resulting mathematical model is given by the following system of differential equations
\begin{eqnarray}
\f{ds}{dt}&=&D(s_{in}-s)-\f{ms}{a+s}x_1-\f{ms}{a+s}x_2+b\nu x_1,\label{1}
\\[1.3ex]
\f{dx_1}{dt}&=&-(\nu+D)x_1+\f{cs}{a+s}x_1-r_1x_1+r_2x_2,\label{2}
\\[1.3ex]
\f{dx_2}{dt}&=&-\nu x_2+\f{cs}{a+s}x_2+r_1x_1-r_2x_2,\label{3}
\end{eqnarray}
\noindent where $s$, $x_1$, $x_2$ denote concentration of the nutrient, the microorganisms in the liquid media and the ones sticked on the walls of the culture vessel, respectively. In addition, $b\in(0,1)$ describes the fraction of dead biomass which is recycled, $\nu>0$ is the collective death rate coefficient of the microbial biomass representing all the aforementioned factors such as diseases, aging, etc. Apart from that, $0<c\leq m$ is the growth rate coefficient of the consumer species. Finally, $r_1>0$ and $r_2>0$ represent the rates at which the species stick on to and shear off from the walls of the culture vessel, respectively.\n

Let us now explain some details about the chemostat model with wall growth. Concerning the equations describing the dynamics of the nutrient or substrate, \eqref{1}, it is quite logical that we need here two terms concerning the consumption of the two different species. In addition, the term $b\nu x_1$ describes the microbial biomass which is recycled and hence can be considered as substrate. In \eqref{2}, the term $-\nu x_1$ reflects the quantity of species which dies, $-D x_1$ denotes the concentration of microbial biomass which is removed from the culture vessel to the collection vessel and, eventually, the last two terms just refer the quantity of microorganisms which stick on to and shear off the walls of the culture vessel. Finally, regarding \eqref{3}, the only term which deserves to be mentioned, since the rest can be explained similarly as before, is the first one. In this case, we can observe that species are sticked on the walls of the culture vessel so they cannot be removed to the collection vessel hence we just write $-\nu x_2$ to take into consideration the microbial biomass death as a consequence of the biological process.\n

There are many works in the literature concerning the deterministic chemostat models (see \cite{caraballo-book,caraballo4,CHKR,MR,HLRS17,J74,HS}). However, some restrictions need to be assumed in these cases, for instance, every parameter is assumed either to be constant or given by a deterministic function depending on the time in the non-autonomous case and these are very strong restrictions, specially if we think that the real world is in most cases random. In this way, our goal is to consider randomness or stochasticity in the chemostat model with wall growth in order to obtain much more realistic mathematical models and get useful information about the long-time behavior of the systems.\n

Nevertheless, there are many different ways of modeling stochasticity or randomness in deterministic models. Concerning the chemostat model, Caraballo {\it et al} \cite{CGL,corrigendumchapter} analyze the simplest chemostat model without wall growth by perturbing the input flow by means of the standard Wiener process. As a result, many drawbacks are found, for instance, it is not possible to prove the positiveness of solution of the resulting stochastic chemostat model due to the huge fluctuations (positive or negative) that the white noise could have. In fact, some state variables could take negative values which is completely unrealistic from the point of view of the applications. In addition, it is not possible to ensure the persistence or coexistence of the species in study either, the main goal pursued by biologists, due to the same reason concerning the white noise.\n

In order to save the previous drawbacks, and motivated by the paper of Imhof and Walcher (see \cite{imhof-walcher}), Caraballo {\it et al} \cite{CGLii} analyze both chemostat models without and with wall growth by means of the standard Wiener process but, in this case, the way of perturbing the model is quite different. More precisely, the authors directly add the noise in every equation of the differential system such that the positiveness of solution could be easily proved (see \cite{CGLii} for details). Nevertheless, since the standard Wiener process has not bounded variation paths, it is not possible to prove the persistence or coexistence of the species in this case either, which is, needless to say, an important drawback also in this case.\n

Thus, Caraballo {\it et al} \cite{CGLR} consider a new noise based on the well-known Ornstein-Uhlenbeck (O-U) process depending on some parameters such that every realization of the noise can be ensured to be bounded inside a positive interval. Thanks to this relevant idea, the authors are able to prove every solution of the simplest chemostat model without wall growth to be positive and, what is much more interesting, they are able to prove the persistence of the total microbial biomass by proving the existence of absorbing and attracting sets which are deterministic, i.e., they do not depend on the realization of the noise. In addition, the dynamics of the resulting random chemostat is analyzed forwards in time, differently to the case when considering the standard Wiener process where a more complicated framework based on the theory of random dynamical systems and pullback attractors is needed.\n

In view of the advantages of using the new suitable Ornstein-Uhlenbeck process presented in \cite{CGLR}, in this paper our aim is to use it to perturb the input flow in the chemostat model with wall growth \eqref{1}-\eqref{3} in the same way that in \cite{CGLR}. As it will be explained after, this allows us to prove the positiveness of solution and we will be able to ensure the persistence of the total microbial biomass and also the coexistence of both species individually, under some conditions on the parameters of the model, in the sense that there exists a number $\eta>0$ such that for any non null initial biomass $x_1(0)$ and $x_2(0)$ each realization satisfies
\begin{equation}
\liminf_{t\rightarrow+\infty} x_i(t)\geq \eta > 0,\label{persistence}
\end{equation}
\noindent for $i=1, 2$, which is stronger than the one used in other many papers $$\liminf_{t\rightarrow +\infty} x(t)>0$$ (see, for instance, \cite{imhof-walcher}).\n

Apart from the previous advantages, we also achieve some improvements when comparing this approach with the one by Xu et al in \cite{XYZ13} since, even though they consider the chemostat model without wall growth, they need the amplitude of the noise to be small enough in order to ensure the persistence of the microbial biomass (see Theorem 1.2 and Section 4) whereas, in our case, modeling the disturbances by using the O-U process, there is no discussion needed on the amplitude of the noise to ensure the persistence (which is, in addition, in the stronger sense \eqref{persistence}). On the other hand, the authors in \cite{XYZ13} prove the results in probability while we prove every results for every realization of the noise in a set of events of full measure.\n

Finally, we will also compare the results obtained with the ones when perturbing the input flow in the chemostat model with wall growth by using the standard Wiener process in order to be able to make a comparison between both ways of modeling.\n

The paper is organized as follows: in Section \ref{s2} we introduce the suitable Ornstein-Uhlenbeck process which was presented in \cite{CGLR}. In Section \ref{s3} we analyze the chemostat model with wall growth \eqref{1}-\eqref{3} where the input flow is perturbed by means of the suitable O-U process. To this end, we prove the existence and uniqueness of global solution which remains inside the positive octant for every positive initial value as well as the existence of a deterministic attracting set whose internal structure will provide us useful information about the long-time behavior of species. Finally, we show some numerical simulations to support the result from the theory. Therefore, in Section \ref{s4} we briefly present the chemostat model with wall growth where the input flow is perturbed by the standard Wiener process and we also display some numerical simulations. Then, in Section \ref{s5} we compare both ways of modeling analyzed in the previous sections in order to see the difference and we also display some numerical simulations. Eventually, in Section \ref{s6} we state some final comments as conclusions of this work.


\section{The Ornstein-Uhlenbeck process}\label{s2}

In this section, we present the key in our current work which simply consists of perturbing the deterministic models by means of a suitable O-U process defined as the following random variable
\begin{equation}
z_{\beta,\gamma}^*(\theta _t\omega )=-\beta\gamma\int\limits_{-\infty }^0e^{\beta s}\theta _t\omega (s)ds,\quad t\in \ensuremath{\mathbb{R}},\,\, \omega \in \Omega,\,\, \beta,\,\,\gamma>0,\label{OU}
\end{equation}
\noindent where $\omega$ denotes a standard Wiener process in a certain probability space $(\Omega,\mathcal{F},\mathbb{P})$, $\beta$ and $\gamma$ are positive parameters which will be explained in more detail below and $\theta_t$ denotes the usual Wiener shift flow given by 
\begin{equation*}
\theta_t \omega(\cdot) = \omega(\cdot + t) - \omega(t),\quad t\in \ensuremath{\mathbb{R}}.\label{ws}
\end{equation*}

We note that the O-U process \eqref{OU} can be obtained as the stationary solution of the Langevin equation
\begin{equation}\label{eq:L}
dz +\beta z dt = \gamma d\omega.
\end{equation}

The O-U process given by \eqref{OU} is a stationary mean-reverting Gaussian stochastic process where $\beta>0$ is a {\it mean reversion constant} that represents the strength with which the process is attracted by the mean or, in other words, how {\it strongly} our system reacts under some perturbation, and $\gamma>0$ is a {\it volatility constant} which represents the variation or the size of the noise independently of the amount of the noise $\alpha>0$. In fact, the O-U process can describe the position of some particle by taking into account the friction, which is the main difference with the standard Wiener process and makes our perturbations to be a better approach to the real ones than the ones obtained when using simply the standard Wiener process. In addition, the O-U process could be understood as a generalization of the standard Wiener process as well in the sense that it would correspond to take $\beta=0$ and $\gamma=1$ in \eqref{OU}. In fact, the O-U also provides a link between the standard Wiener process and no noise at all, as we will see later.\n

Now, we would like to show the relevant effects caused by both parameters $\beta$ and $\gamma$ on the evolution of realization of the O-U process.\n

\begin{figure}[H]
\begin{center}
\includegraphics[width=1.0\textwidth]{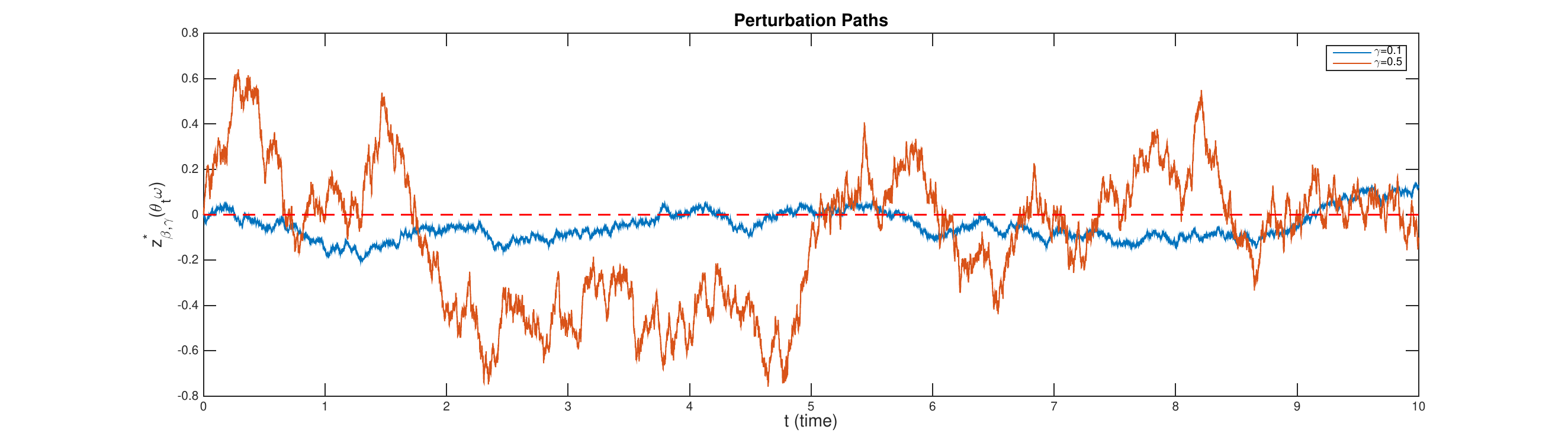}
\caption{Effects of the mean reverting constant on the O-U process}
\label{CGLR-p1}
\end{center}
\end{figure}

{\bf Fixed $\beta>0$.} Then, the volatility of the process increases when considering larger values of $\gamma$ and the evolution of the process is smoother when taking smaller values of $\gamma$, which is completely reasonable due to the fact that $\gamma$ decides the amount of noise introduced to $dz$, the term which measures the variation of the process. Henceforth, the process will be subject to suffer much more disturbances when taking a larger value of $\gamma$. This behavior can be observed in Figure \ref{CGLR-p1}, where we simulate two realizations of the O-U process with $\beta=1$ and we consider $\gamma=0.1$ (blue) and $\gamma=0.5$ (orange).

\begin{figure}[H]
\begin{center}
\includegraphics[width=1.0\textwidth]{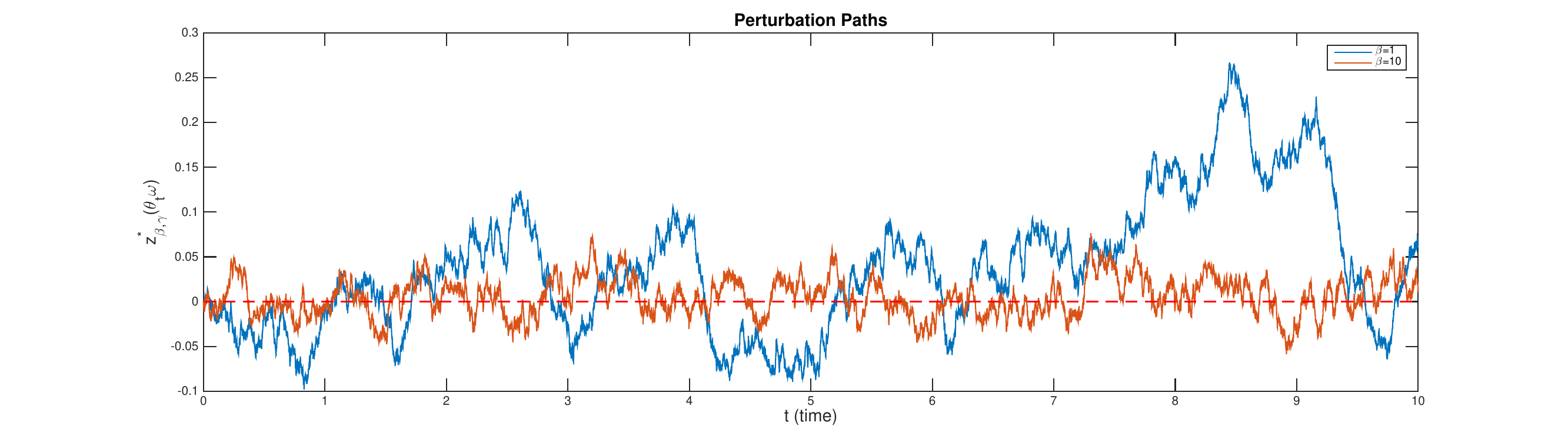}
\caption{Effects of the volatility constant on the O-U process}
\label{CGLR-p2}
\end{center}
\end{figure}

{\bf Fixed $\gamma>0$.} In this case the process tends to go further away from the mean value when considering smaller values of $\beta$ and the attraction of the mean value increases when taking larger values of $\beta$. This behavior is completely logical since $\beta$ has a huge influence on the drift of the Langevin equation \eqref{OU}. We can observe this behavior in Figure \ref{CGLR-p2}, where we simulate two realizations of the O-U process with $\gamma=0.1$ and we take $\beta=1$ (blue) and $\beta=10$ (orange).

Once presented the O-U process and the effects that its parameters cause on its behavior, we state now some essential properties which will be crucial for our analysis in this paper.

\begin{proposition}\label{pOU} 
There exists a $\theta _t$-invariant set $\widetilde{\Omega }\in \mathcal{F}$  of  $\Omega$  of full $\ensuremath{\mathbb{P}}-$measure such that for $\omega \in \widetilde{\Omega }$ and $\beta,\gamma>0$,  we have
\begin{itemize}
\item [(i)] the random variable  $|z_{\beta,\gamma}^*(\omega )|$  is tempered.\n
\item [(ii)] the mapping 
\[
(t,\omega )\rightarrow z_{\beta,\gamma}^*(\theta _t\omega )=-\beta\gamma\int\limits_{-\infty
}^0e^{\beta s}\omega (t+s)\mathrm{d}s+\omega(t)
\]
is a stationary solution of \eqref{OU}
with continuous trajectories;\n
\item [(iii)] for any $\omega \in \tilde \Omega$ one has
\begin{eqnarray*}
\lim_{t\rightarrow \pm \infty }\frac{|z_{\beta,\gamma}^*(\theta _t\omega )|}%
t&=& 0;\\
\lim_{t\rightarrow \pm \infty }\frac 1t\int_0^tz_{\beta,\gamma}^*(\theta _s\omega
)ds&=&0;\\
 \lim_{t\rightarrow \pm \infty }\frac 1t\int_0^t |z_{\beta,\gamma}^*(\theta _s\omega
)| ds&=& \mathbb{E}[z_{\beta,\gamma}^*] < \infty;\n
\end{eqnarray*}
\item [(iv)] finally, for any $\omega\in\widetilde{\Omega}$,
$$\lim_{\beta\rightarrow\infty}z^*_{\beta,\gamma}(\theta_t\omega)=0, \quad\text{for all}\,\, t\in\R.$$
\end{itemize}
\end{proposition} 

The proof of Proposition \ref{pOU} is omitted here. We refer the readers to \cite{alazzawi2017} (Lemma 4.1) for the proof of the last statement and to \cite{arnold,caraballo1} for the proof of the rest.\n

\section{Random chemostat model with wall growth}\label{s3}

In this section, we are interested in analyzing the chemostat model with wall growth \eqref{1}-\eqref{3} where the input flow is perturbed by means of the Ornstein-Uhlenbeck process presented in Section \ref{s2}. We are going to prove the existence and uniqueness of positive global solution of the resulting system as well as the existence of a deterministic attracting set whose internal structure will be also analyzed to obtain detailed information about the long-time behavior of the model. In addition, some conditions on the parameters involved in the system will be provided to ensure the coexistence of the species in the culture vessel.\n

We would also like to note that, for some positive interval, namely $(b_1,b_2)$, thanks to the last property in Proposition \ref{pOU}, for every fixed event $\omega\in\Omega$, it is possible to take $\beta>0$ large enough such that $D+\alpha z^*_{\beta,\delta}(\theta_t\omega)\in(b_1,b_2)$ for every $t\in\R$.

\subsection{Existence and uniqueness of global solution.}

We are interested in analyzing the following random differential system with wall growth and Monod kinetics
\begin{eqnarray}
\f{ds}{dt}&=&\left(D+\azbn\right)(s_{in}-s)-\f{ms}{a+s}x_1-\f{ms}{a+s}x_2+b\nu x_1,\label{rp1}
\\[1.3ex]
\f{dx_1}{dt}&=&-\left(\nu+D+\azbn\right)x_1+\f{cs}{a+s}x_1-r_1x_1+r_2x_2,\label{rp2}
\\[1.3ex]
\f{dx_2}{dt}&=&-\nu x_2+\f{cs}{a+s}x_2+r_1x_1-r_2x_2,\label{rp3}
\end{eqnarray}
\noindent where $z^*_{\beta,\nu}(\theta_t\omega)$ denotes the Ornstein-Uhlenbeck process defined in Section \ref{s2} and every parameter and state variable was introduced in Section \ref{intro}.\n

In this section, $\mathcal{X}=\{(x,y,z)\in\R^3\,:\,x\geq 0,\,y\geq 0, z\geq 0\}$ will denote the positive cone in the three-dimensional space.\n

Firstly, we will state a result concerning the existence and uniqueness of global solution of the chemostat model with wall growth \eqref{rp1}-\eqref{rp3}.

\begin{theorem}
For any initial triple $v_0:=(s_0,x_{10},x_{20})\in\mathcal{X}$, system \eqref{rp1}-\eqref{rp3} possesses a unique global solution $$v(\cdot;0,\omega,v_0):=(s(\cdot;0,\omega,v_0),x_1(\cdot;0,\omega,v_0),x_2(\cdot;0,\omega,v_0))\in\mathcal{C}^1([0,+\infty),\mathcal{X})$$
\noindent with $v(0;0,\omega,v_0)=v_0$, where $s_0:=s(0;0,\omega,v_0)$, $x_{10}:=x_1(0;0,\omega,v_0)$ and $x_{20}:=x_2(0;0,\omega,v_0)$. 
\end{theorem}

\begin{proof}Let us recall that the random system \eqref{rp1}-\eqref{rp3} can be rewritten as
$$\f{dv}{dt}=L(\theta_t\omega)\,\, v+F(v,\theta_t\omega),$$
\noindent where
\begin{eqnarray*}
L(\theta_t\omega) &=& \left(\begin{array}{ccc}
					   -\left(D+\azbn\right) & -m+b\nu & -m \\
					   0 & -\left(\nu+D+\azbn\right)-r_1+c & r_2 \\
					   0 & r_1 & -\nu+c-r_2
					   \end{array}\right)
\end{eqnarray*}
\noindent and $F:\mathcal X\times[0,+\infty)\longrightarrow\R^3$ is given by
\begin{eqnarray*}
F(\eta,\theta_t\omega) &=& \left(\begin{array}{c}
					   \left(D+\azbn\right)s_{in}+\displaystyle{\f{ma}{a+\eta_1}}\eta_2+\displaystyle{\f{ma}{a+\eta_1}}\eta_3\\ \vspace{0.2cm}
					   -\displaystyle{\f{ca}{a+\eta_1}}\eta_2\\ \vspace{0.2cm}
					   -\displaystyle{\f{ca}{a+\eta_1}}\eta_3
					   \end{array}\right),
\end{eqnarray*}
\noindent where $\eta=(\eta_1,\eta_2,\eta_3)\in\mathcal X$.\n

Since $z^*_{\beta,\delta}(\theta_t\omega)$ is continuous, $L$ generates an evolution system on $\R^3$. Moreover, we notice that $F(\cdot,\theta_t\omega)\in\C^1\left(\mathcal X \times[0,+\infty);\R^3\right)$ whence it is locally Lipschitz with respect to $(\eta_1,\eta_2,\eta_3)\in\mathcal{X}$. Thus, system \eqref{rp1}-\eqref{rp3} po\-ssesses a unique local solution.\n

Now, we prove that the unique local solution of system \eqref{rp1}-\eqref{rp3} is defined for any forward time and is, then, a unique global one. To this end, we define the new state variable $p(t)=s(t)+\f{m}{c}(x_1(t)+x_2(t))$ and take into account that $D+\azbn>b_1>0$ for every $t\in\R$, $c\leq m$ and $b\leq 1$. Thus we have that $p$ satisfies the next random differential inequalities
\begin{eqnarray}
\nonumber
\f{dp}{dt}&\leq&\left(D+\azbn\right)s_{in}-b_1s-\f{m}{c}b_1x_1-\f{m}{c}\nu x_2
\\[1.3ex]
\nonumber
&\leq&\left(D+\azbn\right)s_{in}-\vartheta\left[s+\f{m}{c}x_1+\f{m}{c}x_2\right]
\\[1.3ex]
\nonumber
&=&\left(D+\azbn\right)s_{in}-\vartheta p(t),
\end{eqnarray}
\noindent or, in other words, $p$ verifies the following random differential equation
\begin{equation}
\f{dp}{dt}\leq \left(D+\azbn\right)s_{in}-\vartheta p(t),\label{odep}
\end{equation}
\noindent where $\vartheta:=\min\{b_1,\nu\}>0$.\n

By solving \eqref{odep}, we have
\begin{equation}
p(t;0,\omega,p_0)\leq p_0e^{-\vartheta t}+s_{in}\int_0^t\left(D+\azbns\right)e^{-\vartheta(t-s)}ds.\label{solp}
\end{equation}

We remark that the integrand in \eqref{solp} converges to zero for every $t\geq s\geq 0$ when $t$ goes to infinity, but not the integral. Moreover, the integral has subexponential growth.\n

Therefore, $p$ does not blow up at any finite time, thus $s$, $x_1$ and $x_2$ do not blow up at any finite time either. Hence, the solution of system \eqref{rp1}-\eqref{rp3} is defined for any forward time, whence we can straightforwardly deduce that the unique local solution of our system \eqref{rp1}-\eqref{rp3} is, in fact, a unique global one.\n

Now, we are going to prove that the previous unique global solution remains in the positive cone $\mathcal{X}$ for every initial value $v_0\in\mathcal{X}$. To this end, we firstly consider $x_1\geq 0$ and $x_2\geq 0$ and we evaluate the random differential equation for the substrate when $s=0$ such that we have
\begin{equation*}
\left.\f{ds}{dt}\right|_{s=0}=\left(D+\azbn\right)s_{in}+b\nu x_1>0
\end{equation*}
\noindent due to the fact that the perturbed input flow is always positive. Moreover, for every $s\geq 0$ and $x_2\geq 0$, from the equation of the microorganisms in the liquid media, we have
\begin{equation*}
\left.\f{dx_1}{dt}\right|_{x_1=0}=r_2x_2\geq 0
\end{equation*}
\noindent and, for every $s\geq 0$ and $x_1\geq 0$, from the species which are sticked on the walls of the culture vessel, we have
\begin{equation*}
\left.\f{dx_2}{dt}\right|_{x_2=0}=r_1x_1\geq 0.
\end{equation*}

Thus, the unique global solution $v(t;0,\omega,v_0)$ of our random system \eqref{rp1}-\eqref{rp3} remains in the positive cone $\mathcal{X}$ for every initial value $v_0\in\mathcal{X}$.\n

\end{proof}

\subsection{Existence of a deterministic attracting set.}

In this section, we study the existence of a deterministic compact absorbing set as well as the existence of a deterministic attracting set, both of them forwards in time, for the solutions of our random chemostat model with wall growth \eqref{rp1}-\eqref{rp3}. 

\begin{theorem}\label{euap}
For any given $\varepsilon>0$, there exists a deterministic compact absorbing set $B_\varepsilon\subset\mathcal{X}$ for the solutions of system \eqref{rp1}-\eqref{rp3}, i.e., there exists some time $T_F(\omega,\varepsilon)>0$ such that for every given initial pair $v_0\in F$, the solution corresponding to $v_0$ remains inside $B_\varepsilon$ for all $t\geq T_F(\omega,\varepsilon)$.
\end{theorem}

\begin{proof}Consider again the state variable $p(t)=s(t)+\f{m}{c}(x_1(t)+x_2(t))$. Then, from \eqref{solp} we obtain
\begin{eqnarray}
\nonumber
p(t;0,\omega,p_0)&\leq& p_0e^{-\vartheta t}+s_{in}\int_0^t\left(D+\azbns\right)e^{-\vartheta(t-s)}ds
\\[1.3ex]
\nonumber
&\leq& p_0e^{-\vartheta t}+s_{in}\int_0^tb_2e^{-\vartheta(t-s)}ds
\\[1.3ex]
&=& p_0e^{-\vartheta t}+\f{s_{in}b_2}{\vartheta}\left[1-e^{-\vartheta t}\right],\label{labelf}
\end{eqnarray}
\noindent since $D+\alpha z^*_{\beta,\delta}(\theta_s\omega)\leq b_2$ for every $s\in\R$.\n

As a consequence, after making $t$ go to infinity in \eqref{labelf}, we have
\begin{equation}
\lim_{t\rightarrow+\infty}p(t;0,\omega,p_0)\leq \f{s_{in}b_2}{\vartheta}.\label{boundp}
\end{equation}

From \eqref{boundp} we know that, for every initial value $p_0\in F$ and any given $\varepsilon>0$, there exists some time $T_F(\omega,\varepsilon)>0$ such that 
\begin{equation*}
0\leq p(t;0,\omega,p_0)\leq \f{s_{in}b_2}{\vartheta}+\varepsilon
\end{equation*}
\noindent for all $t\geq T_F(\omega,\varepsilon)$. Thus, 
\begin{equation*}
B_\varepsilon=\left\{(s,x_1,x_2)\in\mathcal{X}\,:\, s+\f{m}{c}\left(x_1+x_2\right)\leq \f{s_{in}b_2}{\vartheta}+\varepsilon\right\}\label{absVe}
\end{equation*}
\noindent is, for any $\varepsilon>0$, a deterministic compact absorbing set (forwards in time) for the solutions of system \eqref{rp1}-\eqref{rp3}.\n

\end{proof}

Therefore, thanks to Theorem \ref{euap}, it can be easily deduced that 
\begin{equation}
B_0:=\left\{(s,x_1,x_2)\in\mathcal{X}\,:\, s+\frac{m}{c}\left(x_1+x_2\right)\leq \f{s_{in}b_2}{\vartheta}\right\}\label{b0wg}
\end{equation}
\noindent is a deterministic attracting set (forwards in time) for the solutions of the chemostat model with wall growth \eqref{rp1}-\eqref{rp3}.

\subsection{Internal structure of the deterministic attracting set.}

In this section, our aim is to analyze the internal structure of the deterministic attracting set $B_0$, given by \eqref{b0wg}. To this end, we perform the variable change given by the total biomass and the proportion of the microorganisms in the medium, respectively, as follows
\begin{equation}
x(t)=x_1(t)+x_2(t)\quad\quad\text{and}\quad\quad \xi(t)=\f{x_1(t)}{x_1(t)+x_2(t)}.\label{rxxi}
\end{equation}

For the sake of simplicity, we will write $x$ and $\xi$ instead of $x(t)$ and $\xi(t)$.\n

From \eqref{rxxi}, by differentiation, we obtain the following random differential system
\begin{eqnarray}
\f{ds}{dt}&=&\left(D+\azbn\right)(s_{in}-s)-\f{ms}{a+s}x+b\nu\xi x,\label{rps}
\\[1.3ex]
\f{dx}{dt}&=&-\nu x-\left(D+\azbn\right)\xi x+\f{cs}{a+s}x,\label{rpx}
\\[1.3ex]
\f{d\xi}{dt}&=&-\left(D+\azbn\right)\xi(1-\xi)-r_1\xi+r_2(1-\xi).\label{rpxi}
\end{eqnarray}

We remark that the dynamics of the proportion of the species in the liquid media, $\xi$, is uncoupled of the rest of the system, then we first analyze its asymptotic behavior and we investigate the rest of the system in a second step.\n

Thanks to \eqref{rxxi}, it is straightforward to prove by definition that
\begin{equation*}
0\leq \xi(t;0,\omega,\xi_0)\leq 1\label{bxi}
\end{equation*}
\noindent for every $t\geq 0$ and any initial value $\xi_0\in(0,1)$. In addition, from \eqref{rpxi}, we can evaluate the corresponding random differential equation when $\xi=0$ and $\xi=1$, respectively, such that we obtain
\begin{equation*}
\left.\f{d\xi}{dt}\right|_{\xi=0}=r_2>0\quad\quad\quad\text{and}\quad\quad\quad\left.\f{d\xi}{dt}\right|_{\xi=1}=-r_1<0,
\end{equation*}
\noindent whence we notice that the interval $(0,1)\subset\R$ defines a positively invariant set for the dynamics of the proportion.\n

On the one hand, thanks to the fact that $b_1<D+\azbn<b_2$ for every $t\in\R$, from \eqref{rpxi} we have
\begin{eqnarray}
\nonumber
\f{d\xi}{dt}&=&-\left(D+\azbn\right)\xi(1-\xi)-r_1\xi+r_2(1-\xi)
\\[1.3ex]
\nonumber
&\leq&-(b_1+r_1+r_2)\xi+b_1+r_2,
\end{eqnarray}
\noindent hence we obtain the following random differential equation
\begin{equation}
\f{d\xi}{dt}\leq -(b_1+r_1+r_2)\xi+b_1+r_2.\label{odexiupper}
\end{equation}

By solving now \eqref{odexiupper}, we obtain the following upper bound
\begin{equation*}
\xi(t;0,\omega,\xi_0)\leq\xi_0e^{-(b_1+r_1+r_2)t}+\f{b_1+r_2}{b_1+r_1+r_2}\left[1-e^{-(b_1+r_2+r_2)t}\right]
\end{equation*}
\noindent for any initial value $\xi_0\in(0,1)$ and for all $t\geq 0$.\n

On the other hand, from \eqref{rpxi} we also have
\begin{eqnarray}
\nonumber
\f{d\xi}{dt}&=&-\left(D+\azbn\right)\xi(1-\xi)-r_1\xi+r_2(1-\xi)
\\[1.3ex]
\nonumber
&\geq&-(b_2+r_1+r_2)\xi+r_2,
\end{eqnarray}
\noindent whence we obtain
\begin{equation}
\f{d\xi}{dt}\geq -(b_2+r_1+r_2)\xi+r_2.\label{odexilower}
\end{equation}

By solving in this case \eqref{odexilower}, we obtain the following lower bound
\begin{equation*}
\xi(t;0,\omega,\xi_0)\geq\xi_0e^{-(b_2+r_1+r_2)t}+\f{r_2}{b_2+r_1+r_2}\left[1-e^{-(b_2+r_1+r_2)t}\right]
\end{equation*}
\noindent for any initial value $\xi_0\in(0,1)$ and for all $t\geq 0$.\n

From the calculations above, we have the following bounds for the dynamics of the proportion $\xi$, which are given by
\begin{equation}
\xi(t;0,\omega,\xi_0)\leq\xi_0e^{-(b_1+r_1+r_2)t}+\f{b_1+r_2}{b_1+r_1+r_2}\left[1-e^{-(b_1+r_2+r_2)t}\right]\label{boundsxiu}
\end{equation}
\noindent and
\begin{equation}
\xi(t;0,\omega,\xi_0)\geq\xi_0e^{-(b_2+r_1+r_2)t}+\f{r_2}{b_2+r_1+r_2}\left[1-e^{-(b_2+r_1+r_2)t}\right]\label{boundsxil}
\end{equation}
\noindent for any initial value $\xi_0\in(0,1)$ and for all $t\geq 0$.\n

Then, by making $t$ go to infinity in \eqref{boundsxiu} and \eqref{boundsxil}, respectively, we obtain
\begin{equation*}
\f{r_2}{b_2+r_1+r_2}\leq\lim_{t\rightarrow+\infty}\xi(t;0,\omega,\xi_0)\leq\f{b_1+r_2}{b_1+r_1+r_2}
\end{equation*}
\noindent for every any initial value $\xi_0\in(0,1)$.\n

Hence, for any given $\varepsilon>0$ and any initial value $\xi_0\in(0,1)$, there exists some time $T(\omega,\varepsilon)>0$ such that
\begin{equation*}
-\varepsilon+\f{r_2}{b_2+r_1+r_2}\leq\xi(t;0,\omega,\xi_0)\leq\f{b_1+r_2}{b_1+r_1+r_2}+\varepsilon
\end{equation*}
\noindent for all $t\geq T(\omega,\varepsilon)$.\n

Thus,
\begin{equation*}
B_\varepsilon^\xi=\left\{\xi\in(0,1)\,:\, -\varepsilon+\xi^*_l\leq\xi\leq\xi^*_u+\varepsilon\right\}\label{absxi}
\end{equation*}
\noindent defines a deterministic compact absorbing set for the dynamics of the proportion, where $\xi^*_l$ and $\xi^*_u$ are both deterministic constants given by
\begin{equation}
\xi^*_l:=\f{r_2}{b_2+r_1+r_2}\quad\quad\text{and}\quad\quad \xi^*_u:=\f{b_1+r_2}{b_1+r_1+r_2}.\label{boundsxiat}
\end{equation}

As a consequence, the dynamics of the proportion remains asymptotically inside $B^\xi_\varepsilon$ for any given $\varepsilon>0$ and, then, we obtain the following attracting set for the corresponding state variable describing the dynamics of the proportion
\begin{equation}
B_0^\xi:=\left\{\xi\in(0,1)\,:\,\xi^*_l\leq\xi\leq\xi^*_u\right\}.
\end{equation}

We remark that, since the constants defined in \eqref{boundsxiat} are deterministic, both $B^\xi_\varepsilon$ and $B_0^\xi$ are also deterministic sets, i.e., they do not depend on the noise. In addition, they are absorbing sets forwards in time.\n

Now, we focus on the analysis of the dynamics of the substrate, $s$, and the microorganisms concentration, $x$. We already proved that, for every time $t$ large enough, the dynamics of the proportion satisfies the following inequalities
\begin{equation}
\xi^*_l\leq \xi(t;0,\omega,\xi_0)\leq \xi^*_u\label{bx}
\end{equation}
\noindent for every initial value $\xi_0\in(0,1)$.\n

Having reached this point, we define a new state variable
\begin{equation*}
z(t)=cs(t)+mx(t).\label{rz}
\end{equation*}

We will write $z$, instead of $z(t)$, for the sake of simplicity.\n

Hence, by differentiation, due to the fact that $b\leq 1$, $c\leq m$ and since $\xi(t;0,\omega,\xi_0)\leq 1$ for every $t\geq 0$ and any initial value $\xi_0\in(0,1)$, thanks to \eqref{bx}, we obtain that $z$ satisfies the following random differential equations
\begin{equation}
\f{dz}{dt}\leq-\left(D+\azbn\right)\xi^*_lz+cs_{in}\left(D+\azbn\right)\label{odezupper}
\end{equation}
\noindent and
\begin{equation}
\f{dz}{dt}\geq-\left[\nu+\left(D+\azbn\right)-\f{cb\nu}{m}\xi^*_l\right]z+c\left(D+\azbn\right)s_{in}\label{odezlower}
\end{equation}
\noindent for every time $t$ large enough.\n

By solving now \eqref{odezupper} and \eqref{odezlower},  thanks to bounds of the O-U process, we obtain
\begin{eqnarray}
\nonumber
z(t;0,\omega,z_0)&\leq& z_0e^{-D\xi^*_lt-\alpha\xi^*_l\int_0^tz^*_{\beta,\delta}(\theta_r\omega)dr}
\\[1.3ex]
\nonumber
&&+cs_{in}\int_0^t\left(D+\azbns\right)e^{-D\xi^*_l(t-s)-\alpha\xi^*_l\int_s^tz^*_{\beta,\delta}(\theta_r\omega)dr}ds
\\[1.3ex]
&\leq& z_0e^{-b_1\xi^*_lt}+\f{cs_{in}b_2}{\xi^*_lb_1}\left[1-e^{-\xi^*_lb_1t}\right]\label{boundzupper}
\end{eqnarray}
\noindent and
\begin{eqnarray}
\nonumber
z(t;0,\omega,z_0)&\geq& z_0e^{-\left(\nu+D-\f{cb\nu}{m}\xi^*_l\right)t-\alpha\int_0^tz^*_{\beta,\delta}(\theta_r\omega)dr}
\\[1.3ex]
\nonumber
&&+cs_{in}\int_0^t\left(D+\azbns\right)e^{-\left(\nu+D-\f{cb\nu}{m}\xi^*_l\right)(t-s)-\alpha\int_s^tz^*_{\beta,\delta}(\theta_r\omega)dr}ds
\\[1.3ex]
&\geq& z_0e^{-\left(b_2+\nu-\f{cb\nu}{m}\xi^*_l\right)t}+\f{cs_{in}b_1}{b_2+\nu-\f{cb\nu}{m}\xi^*_l}\left[1-e^{-\left(b_2+\nu-\f{cb\nu}{m}\xi^*_l\right)t}\right],\label{boundzlower}
\end{eqnarray}
\noindent respectively, for every time $t$ large enough.\n


Thus, after making $t$ go to infinity in \eqref{boundzupper} and \eqref{boundzlower}, we have
\begin{equation}
z^*_l:=\f{cs_{in}b_1}{b_2+\nu-\f{cb\nu}{m}\xi^*_l}\leq\lim_{t\rightarrow+\infty}z(t;0,\omega,z_0)\leq \f{cs_{in}b_2}{\xi^*_lb_1}=:z^*_u,\label{ze}
\end{equation}
\noindent for every initial value $v_0\in F$, where we used the fact that $b_2+\nu-\f{cb\nu}{m}\xi^*_l>0$ is fulfilled.\n

We would like to remark that both constants $z^*_l$ and $z^*_u$ in \eqref{ze} do not depend on the noise $\omega$ or, in oder words, we obtained in \eqref{ze} upper and lower deterministic bounds for the dynamics of $z$, what is more, forwards in time.\n

From \eqref{ze}, we have that, for every $v_0\in F$ and any $\varepsilon>0$, there exists some time $T_F(\omega,\epsilon)>0$ such that
\begin{equation}
z^*_l-\varepsilon\leq z(t;0,\omega,z_0)\leq z^*_u+\varepsilon\label{bzc}
\end{equation}
\noindent holds true for all $t\geq T_F(\omega,\varepsilon)$.\n

As a result, we deduce that, for any $\varepsilon>0$,

\begin{equation}
B^{(s,x)}_\varepsilon:=\left\{(s,x)\in\R^2_+\,:\, z^*_l-\varepsilon\leq cs+mx\leq z^*_u+\varepsilon\right\}\label{sasx}
\end{equation}
\noindent is a deterministic compact absorbing set (forwards in time) for the solutions of system \eqref{rps}-\eqref{rpx}.\n

\begin{figure}[H]
\begin{center}
\includegraphics[width=0.8\textwidth]{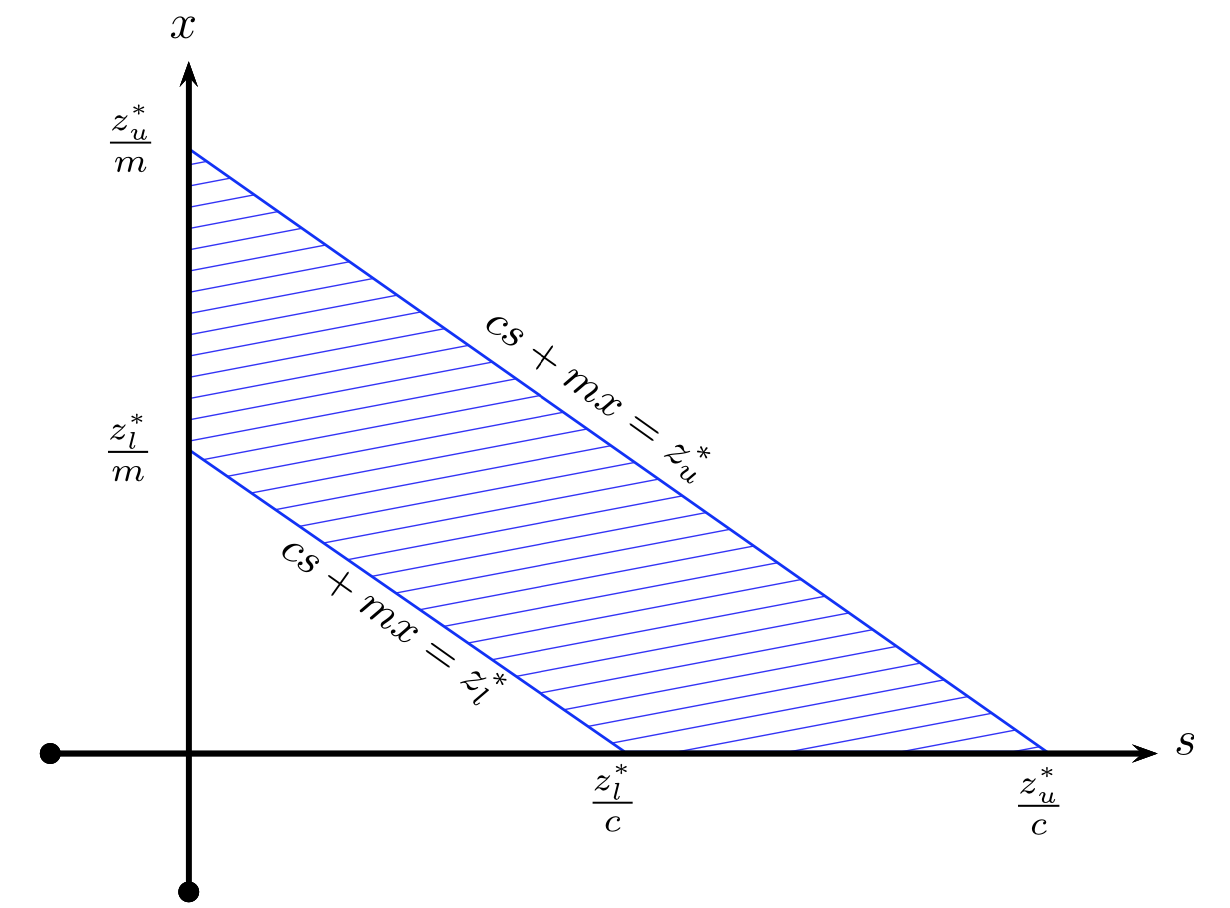}
\end{center}
\caption{Attracting set $B^{(s,x)}_0$}
\label{fabszrbb}
\end{figure}

Therefore, we obtain the following attracting set (forwards in time) for the solutions of system \eqref{rps}-\eqref{rpx}
\begin{equation*}
B^{(s,x)}_0:=\left\{(s,x)\in\R^2_+\,:\, z^*_l\leq cs+mx\leq z^*_u\right\},
\end{equation*}
\noindent see Figure \ref{fabszrbb}.

Now we analyze the dynamics of both the nutrient and the species individually in order to obtain more detailed information about the long-time behavior of the random system \eqref{rps}-\eqref{rpx}. In addition, we will provide conditions under which the persistence in the {\it strong} sense \eqref{persistence} of both species, the ones in the medium and the ones sticked on the walls of the culture vessel, can be proved.\n

\begin{proposition}
Assume that the following condition
\begin{equation}
\nu+D\xi^*_l>c\label{ce}
\end{equation}
\noindent holds true. Then, the attracting set for the solutions of the chemostat model with wall growth \eqref{rp1}-\eqref{rp3} is reduced to a deterministic segment, more precisely, it is
\begin{equation*}
\wh{B}_0^{(s,x)}=\left[\f{z^*_l}{c},\f{z^*_u}{c}\right]\times\{0\}\times\{0\}.
\end{equation*}
\end{proposition}

\begin{proof}On the one hand, from \eqref{rpx} we have that $x$ satisfies the following random differential inequality
\begin{equation*}
\f{dx}{dt}\leq-\left[\nu+\left(D+\azbn\right)\xi^*_l-c\right]x,\label{odexe}
\end{equation*}
\noindent for every time $t$ large enough, whose solution is given by
\begin{equation}
x(t;0,\omega,x_0)\leq x_0e^{-(\nu+D\xi^*_l-c)t-\alpha\xi^*_l\int_0^tz^*_{\beta,\delta}(\theta_s\omega)ds}.\label{boundxe}
\end{equation}

Besides, from \eqref{boundxe}, we have that
\begin{equation*}
\lim_{t\rightarrow+\infty}x(t;0,\omega,x_0)\leq 0
\end{equation*}
\noindent as long as \eqref{ce} is fulfilled or, in other words, both species become extinct if \eqref{ce} holds true.\n

\end{proof}

The next result states a condition for the persistence of the microorganisms in the strong sense \eqref{persistence} can be proved.

\begin{theorem}
Assume that
\begin{equation}
\nu+b_2<\f{z^*_l}{a+\f{z^*_u}{c}}\label{cond}
\end{equation}
\noindent is fulfilled and \eqref{ce} does not hold, where we recall that $z^*_l$ and $z^*_u$ are the constants defined as
\begin{equation*}
z^*_l:=\f{cs_{in}b_1}{b_2+\nu-\f{cb\nu}{m}\xi^*_l}\quad\quad\text{and}\quad\quad z^*_u:=\f{cs_{in}b_2}{\xi^*_lb_1}.\label{defes}
\end{equation*}
\noindent Then, there exists a deterministic compact absorbing set, which is strictly contained in the first quadrant of the two-dimensional space, for the solutions of system \eqref{rps}-\eqref{rpx}.
\end{theorem}

\begin{proof}Let us recall the random differential equation describing the dynamics of the microorganisms concentration, $x=x_1+x_2$, which is given by
\begin{equation*}
\f{dx}{dt}=-\nu x-\left(D+\azbn\right)\xi x+\f{cs}{a+s}x,
\end{equation*}
\noindent whence we obtain that
\begin{eqnarray}
\nonumber
\f{dx(t;0,\omega,x_0)}{dt}&=&-\nu x(t;0,\omega,x_0)-\left(D+\alpha z^*_{\beta,\delta}(\theta_t\omega)\right)\xi(t;0,\omega,\xi_0)x(t;0,\omega,x_0)
\\[1.3ex]
&&+\f{cs(t;0,\omega,s_0)}{a+s(t;0,\omega,s_0)}x(t;0,\omega,x_0).\label{xc}
\end{eqnarray}

On the one hand, thanks to the definition of the proportion \eqref{rxxi}, we have that 
\begin{equation}
0\leq \xi(t;0,\omega,\xi_0)\leq 1\label{bxna}
\end{equation}
\noindent for every $t\geq 0$ and any initial value $\xi_0\in(0,1)$.\n

Thus, from \eqref{xc}, thanks to the previous calculations and the bound of the O-U process, we obtain the following random differential inequality
\begin{eqnarray}
\nonumber
\f{dx(t;0,\omega,x_0)}{dt}&\geq&-\nu x(t;0,\omega,x_0)-b_2x(t;0,\omega,x_0)
\\[1.3ex]
&&+\f{cs(t;0,\omega,s_0)}{a+s(t;0,\omega,s_0)}x(t;0,\omega,x_0)\label{xci}
\end{eqnarray}
\noindent for all $t\geq 0$ and every initial value $v_0\in F$.\n

By definition, we know that
\begin{equation*}
z(t;0,\omega,z_0)=cs(t;0,\omega,s_0)+mx(t;0,\omega,x_0)
\end{equation*}
\noindent and, thanks to \eqref{bzc}, we have that, for each $v_0\in F$ and any $\varepsilon>0$, there exists some time $T_F(\omega,\varepsilon)>0$ such that
\begin{equation*}
z^*_l-\varepsilon\leq cs(t;0,\omega,s_0)+mx(t;0,\omega,x_0)\leq z^*_u+\varepsilon
\end{equation*}
\noindent holds true for every $t\geq T_F(\omega,\varepsilon)$.\n

As a consequence, since $c\leq m$, we have the following inequalities
\begin{equation*}
cs(t;0,\omega,s_0)\geq z^*_l-\varepsilon-mx(t;0,\omega,x_0)\label{slowercd}
\end{equation*}
\noindent and
\begin{equation*}
s(t;0,\omega,s_0)\leq\f{z^*_u}{c}+\f{\varepsilon}{c}-x(t;0,\omega,x_0)
\end{equation*}
\noindent for every initial value $v_0\in F$, any $\varepsilon>0$ and for all $t\geq T_F(\omega,\varepsilon)$.\n

Then, from \eqref{xci}, we have 
\begin{eqnarray}
\nonumber
\f{dx(t;0,\omega,x_0)}{dt}&\geq&-\nu x(t;0,\omega,x_0)-b_2x(t;0,\omega,x_0)
\\[1.3ex]
&&+\f{z^*_l-mx(t;0,\omega,x_0)-\varepsilon}{a+\f{z^*_u}{c}+\f{\varepsilon}{c}-x(t;0,\omega,x_0)}x(t;0,\omega,x_0)\label{odexb}
\end{eqnarray}
\noindent for every $v_0\in F$, any $\varepsilon>0$ and for all $t\geq T_F(\omega,\varepsilon)$.\n

Now, we study the differential equation \eqref{odexb} when $x=\widetilde{x}$, where $\widetilde{x}$ is defined as
\begin{equation}
\widetilde{x}=\f{z^*_l-\left(\nu+b_2\right)\left(a+\f{z^*_u}{c}\right)}{m+c}.\label{xt}
\end{equation}

Then, from \eqref{odexb} and considering $\varepsilon<c\widetilde{x}$, we obtain
\begin{eqnarray}
\nonumber
\left.\f{dx(t;0,\omega,x_0)}{dt}\right|_{x=\widetilde{x}}&\geq&\left[-(\nu+b_2)+\f{z^*_l-m\widetilde{x}-\varepsilon}{a+\f{z^*_u}{c}+\f{\varepsilon}{c}-\widetilde{x}}\right]\widetilde{x}
\\[1.3ex]
\nonumber
&>&\left[-(\nu+b_2)+\f{z^*_l-m\widetilde{x}-c\widetilde{x}}{a+\f{z^*_u}{c}+\f{c\widetilde{x}}{c}-\widetilde{x}}\right]\widetilde{x}=0
\end{eqnarray}
\noindent for every $v_0\in F$, any $\varepsilon\in(0,c\widetilde{x})$ and for all $t\geq T_F(\omega,\varepsilon)$.\n

Hence, as long as \eqref{cond} is fulfilled, we have that, for any $\varepsilon\in(0,c\widetilde{x})$, where $\widetilde{x}$ is given by \eqref{xt}, for every $v_0\in F$, there exists some time $T_F(\omega,\varepsilon)>0$ such that
\begin{eqnarray}
\left.\f{dx(t;0,\omega,x_0)}{dt}\right|_{x=\widetilde{x}}&>&0\label{xtb}
\end{eqnarray}
\noindent for all $t\geq T_F(\omega,\varepsilon)$.\n

Therefore, from \eqref{xtb} we conclude that, as long as \eqref{cond} is fulfilled, we have the following lower deterministic bound for the dynamics of the species
\begin{equation*}
x(t;0,\omega,x_0)>\widetilde{x},\label{xboundf}
\end{equation*}
\noindent for any $\varepsilon\in(0,c\widetilde{x})$, every given $v_0\in F$ and for all $t\geq T_F(\omega,\varepsilon)$.\n

Now, let us recall the random differential equation held by the substrate
\begin{equation*}
\f{ds}{dt}=\left(D+\azbn\right)(s_{in}-s)-\f{ms}{a+s}x+b\nu\xi x
\end{equation*}
\noindent for every $t\geq 0$, whence we obtain that
\begin{eqnarray}
\nonumber
\f{ds(t;0,\omega,s_0)}{dt}&=&\left(D+\alpha z^*_{\beta,\delta}(\theta_t\omega)\right)s_{in}-\left(D+\alpha z^*_{\beta,\delta}(\theta_t\omega)\right)s(t;0,\omega,s_0)
\\[1.3ex]
\nonumber
&&-\f{ms(t;0,\omega,s_0)}{a+s(t;0,\omega,s_0)}x(t;0,\omega,x_0)+b\nu\xi(t;0,\omega,\xi_0)x(t;0,\omega,x_0),\label{sc}
\end{eqnarray}
\noindent for all $t\geq 0$ and every initial value $v_0\in F$.\n

Moreover, from \eqref{bzc} and \eqref{bxna}, we know that, for each $v_0\in F$ and any $\varepsilon>0$, there exists some time $T_F(\omega,\varepsilon)>0$ such that
\begin{equation*}
x(t;0,\omega,x_0)\leq \f{z^*_u}{m}+\f{\varepsilon}{m}\label{xuppercd}
\end{equation*}
\noindent holds true for every $t\geq T_F(\omega,\varepsilon)$.\n

Summing up, thanks to the previous calculations, from \eqref{sc} we have 
\begin{eqnarray}
\nonumber
\f{ds(t;0,\omega,s_0)}{dt}&\geq&b_1s_{in}-b_2s(t;0,\omega,s_0)-\f{ms(t;0,\omega,s_0)}{a+s(t;0,\omega,s_0)}\f{z^*_u+\varepsilon}{m},\label{sci}
\end{eqnarray}
\noindent for every $v_0\in F$, any $\varepsilon>0$ and for all $t\geq T_F(\omega,\varepsilon)$.\n

Now, we study the differential equation \eqref{sci} when $s=\widetilde{s}$, where $\widetilde{s}$ is defined as
\begin{equation}
\widetilde{s}=\f{b_1s_{in}}{b_2+2\f{z^*_u}{a}}.\label{stbb}
\end{equation}

Then, from \eqref{sci} and considering $\varepsilon<z^*_u$, we obtain
\begin{eqnarray}
\nonumber
\left.\f{ds(t;0,\omega,s_0)}{dt}\right|_{s=\widetilde{s}}&\geq&b_1s_{in}-b_2\widetilde{s}-\f{\widetilde{s}}{a+\widetilde{s}}(z^*_u+\varepsilon)
\\[1.3ex]
\nonumber
&>&b_1s_{in}-b_2\widetilde{s}-\f{2\widetilde{s}}{a}z^*_u=0
\end{eqnarray}
\noindent for every $v_0\in F$, any $\varepsilon\in(0,z^*_u)$ and for all $t\geq T_F(\omega,\varepsilon)$.\n

Hence, we have that, for any $\varepsilon\in(0,z^*_u)$ and every $v_0\in F$, there exists some time $T_F(\omega,\varepsilon)>0$ such that
\begin{eqnarray*}
\left.\f{ds(t;0,\omega,s_0)}{dt}\right|_{s=\widetilde{s}}&>&0
\end{eqnarray*}
\noindent for all $t\geq T_F(\omega,\varepsilon)$.\n

Thus, we obtain the following lower deterministic bound for the dynamics of the substrate
\begin{equation*}
s(t;0,\omega,s_0)>\widetilde{s},\label{sboundf}
\end{equation*}
\noindent for any $\varepsilon\in(0,z^*_u)$, every given $v_0\in F$ and for all $t\geq T_F(\omega,\varepsilon)$.\n

In conclusion, we obtain that, by taking any $\varepsilon\in(0,\min\{c\widetilde{x},z^*_u\})$, where $\widetilde{x}$ is given by \eqref{xt}, for every given $v_0\in F$, there exists some time $T_F(\omega,\varepsilon)>0$ such that
\begin{equation}
x(t;0,\omega,x_0)>\widetilde{x}\label{bxf}
\end{equation}
\noindent and
\begin{equation}
s(t;0,\omega,s_0)>\widetilde{s}\label{bxfb}
\end{equation}
\noindent hold true for all $t\geq T_F(\omega,\varepsilon)$.\n

As a result, we can deduce that, for any $\varepsilon\in(0,\min\{c\widetilde{x},z^*_u\})$,
\begin{equation}
\wh{B}^{(s,x)}_\varepsilon=\left\{(s,x)\in\R^2_+\,:\, x\geq\widetilde{x},\,\, s\geq\widetilde{s},\,\, z^*_l-\varepsilon\leq cs+mx\leq z^*_u+\varepsilon\right\},\label{sxso}
\end{equation}
\noindent where $\widetilde{x}$ and $\widetilde{s}$ are defined by \eqref{xt} and \eqref{stbb}, respectively, is a deterministic compact absorbing set (forwards in time) for the solutions of system \eqref{rps}-\eqref{rpx}.\n

\end{proof}

It is worth mentioning that we already proved the dynamics of system \eqref{rps}-\eqref{rpx} to remain inside $B^{(s,x)}_\varepsilon$, defined as in \eqref{sasx}, forwards in time. In the previous result, as long as \eqref{cond} is fulfilled, we obtain in addition a smaller deterministic compact absorbing set $\widetilde{B}^{(s,x)}_\varepsilon$ forwards in time as well, defined by \eqref{sxso}, which besides is strictly contained in the first quadrant of the two-dimensional space. This fact will be the main key to guarantee the persistence of both species, individually, in the strong sense \eqref{persistence}.\n

\begin{figure}[H]
\begin{center}
\includegraphics[width=0.8\textwidth]{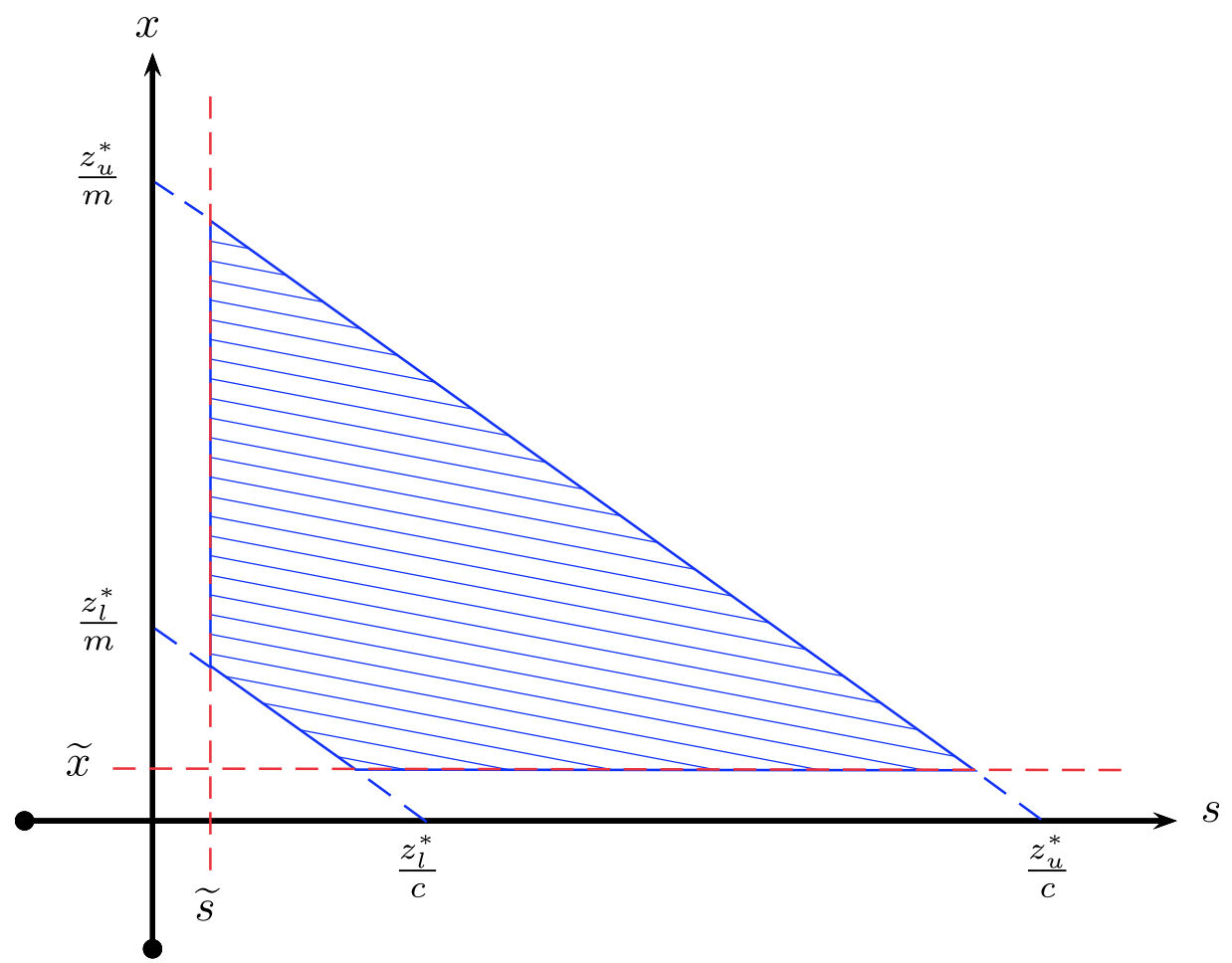}
\end{center}
\caption{Attracting set $\widetilde{B}^{(s,x)}_0$}
\label{abstb}
\end{figure}

Therefore,
\begin{equation*}
\widetilde{B}^{(s,x)}_0:=\left\{(s,x)\in\R^2_+\,\,:\,\,x\geq\widetilde{x},\,\,s\geq\widetilde{s},\,\,z^*_l\leq cs+mx\leq z^*_u\right\}
\end{equation*}
\noindent is a deterministic attracting set (forwards in time) for the solutions of system \eqref{rps}-\eqref{rpx}, see Figure \ref{abstb}.

\begin{remark} 
It is not difficult to check that both $\widetilde{x}<\f{z^*_l}{m}$ and $\widetilde{s}<\f{z^*_l}{c}$ are satisfied.
\end{remark}

Finally, we will analyze the dynamics of both species, $x_1$ and $x_2$, individually, to prove that both of them also persist as long as \eqref{cond} holds true. To this end, thanks to \eqref{bxf}, \eqref{bxfb} and the definition of the proportion, $\xi=x_1/x$, we obtain that
\begin{eqnarray}
\nonumber
x_1(t;0,\omega,x_{10})&=&\xi(t;0,\omega,\xi_0)x(t;0,\omega,x_0)
\\[1.3ex]
\nonumber
&>&\xi^*_l\widetilde{x}>0\label{x1f}
\end{eqnarray}
\noindent for every $t$ large enough and any initial value $v_0\in F$.\n

In addition, we also have
\begin{eqnarray}
\nonumber
x_2(t;0,\omega,x_{20})&=&x(t;0,\omega,x_0)(1-\xi(t;0,\omega,\xi_0))
\\[1.3ex]
\nonumber
&>&(1-\xi^*_u)\widetilde{x}>0\label{x2f}
\end{eqnarray}
\noindent for every $t$ large enough and any initial value $v_0\in F$.\n

Hence, we obtain that both species, the ones in the medium and also the ones sticked on to the walls of the culture vessel, will persist as long as \eqref{cond} holds true.

\begin{remark}
It is possible to improve the deterministic lower bounds obtained in \eqref{bxf} and \eqref{bxfb} by considering smaller values of $\varepsilon>0$. Particularly, we could consider $\varepsilon\in\left(0,\f{c}{n}\widetilde{x}\right)$, for any $n\in\mathbb{N}$, instead of $\varepsilon\in(0,c\widetilde{x})$ such that we have that $x>\widetilde{x}_n$, instead of $x>\widetilde{x}$ as in \eqref{bxf}, where $\widetilde{x}_n$ is given by
\begin{equation*}
\widetilde{x}_n:=\f{z^*_l-(\nu+b_2)\left(a+\f{z^*_u}{c}\right)}{m+\f{c}{n}},
\end{equation*}
\noindent which clearly satisfies
\begin{equation*}
\widetilde{x}_n:=\f{z^*_l-(\nu+b_2)\left(a+\f{z^*_u}{c}\right)}{m+\f{c}{n}}>\f{z^*_l-(\nu+b_2)\left(a+\f{z^*_u}{c}\right)}{m+c}=:\widetilde{x}>0,
\end{equation*}
\noindent for any $n\in\mathbb{N}$, since $n\geq 1$.\n

Similarly, we could consider $\varepsilon\in\left(0,\f{1}{n}z^*_u\right)$, for any $n\in\mathbb{N}$, instead of $\varepsilon\in(0,z^*_u)$ such that we get $s>\widetilde{s}_n$ instead of $s>\widetilde{s}$ as in \eqref{bxfb}, where $\widetilde{s}_n$ is given by
\begin{equation*}
\widetilde{s}_n:=\f{b_1s_{in}}{b_2+\f{z^*_u}{a}\left(1+\f{1}{n}\right)}
\end{equation*}
\noindent which clearly verifies
\begin{equation*}
\widetilde{s}_n:=\f{b_1s_{in}}{b_2+\f{z^*_u}{a}\left(1+\f{1}{n}\right)}>\f{b_1s_{in}}{b_2+\f{2z^*_u}{a}}=:\widetilde{s}>0,
\end{equation*}
\noindent for any $n\in\mathbb{N}$, since $n\geq 1$.
\end{remark}

\subsection{Numerical simulations.}

In this section we show several numerical simulations for different values of the parameters involved in the random chemostat model with wall growth \eqref{rp1}-\eqref{rp3} in order to support the results provided through this section. The blue dashed lines represent the solution of the corresponding deterministic chemostat model whereas the other ones represent different realizations of the solution of the random system. Moreover, we display four different panels in each figure: there is a big one on the left-hand side showing the general dynamics of the model and there are three smaller panels on the right-hand side where the individual dynamics of the substrate and both species will be presented.\n 

\begin{figure}[H]
\begin{center}
\includegraphics[scale=0.27]{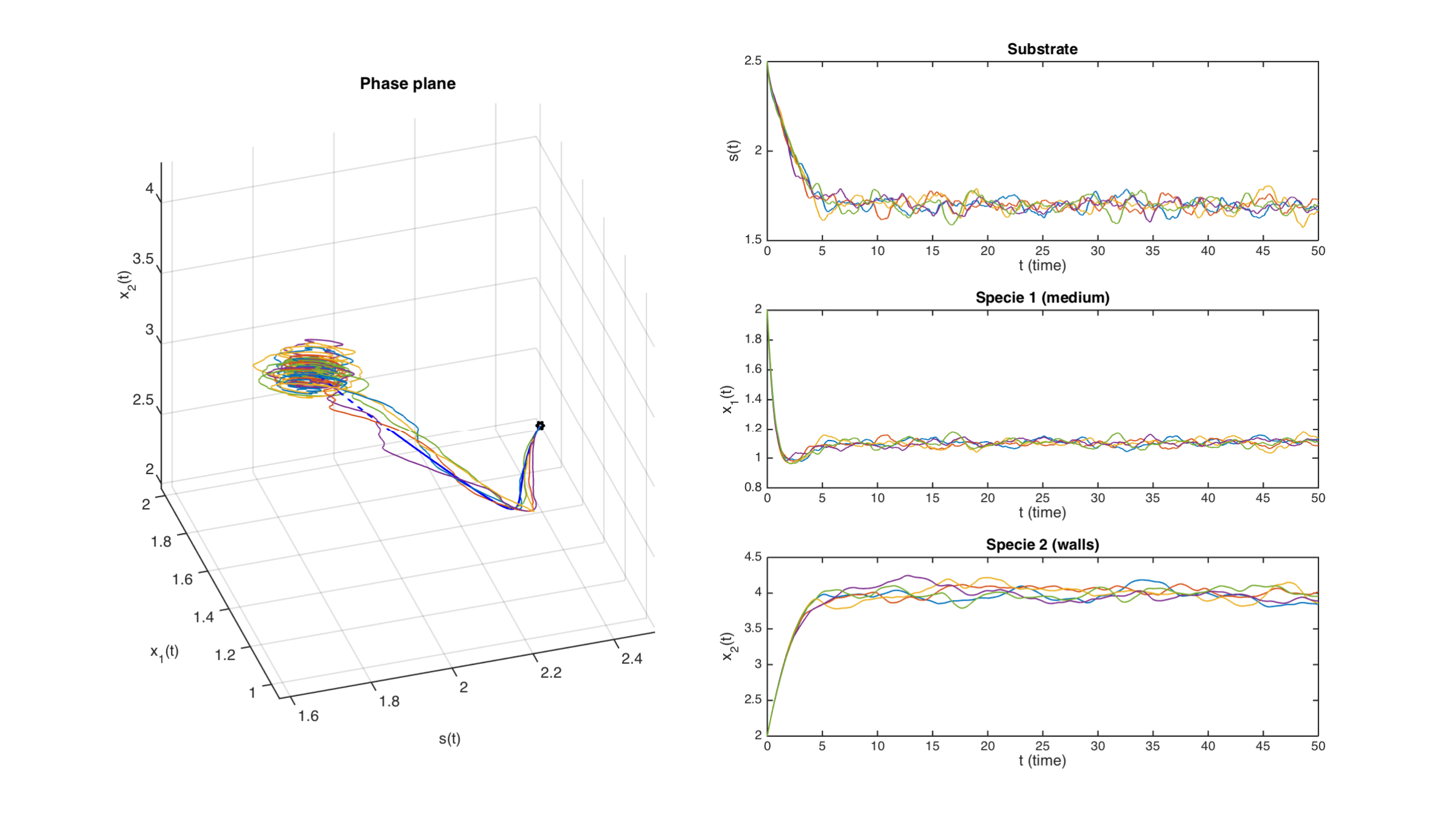}
\end{center}
\caption{Persistence of the species in the random chemostat model with wall growth}
\label{2_11}
\end{figure}

On the one hand, we consider in Figure \ref{2_11} $s_{in}=4$, $D=2$, $a=1.6$, $m=2$, $b=0.5$, $\nu=1.2$, $c=3$, $r_1=0.2$, $r_2=0.4$, $\alpha=0.5$, $\beta=1$, $\nu=0.2$ and we will take $s_0=2.5$, $x_{10}=2$, $x_{20}=2$ as initial values for the substrate and both species, respectively. As a result, we can see that both species persist.

\begin{figure}[H]
\begin{center}
\includegraphics[scale=0.27]{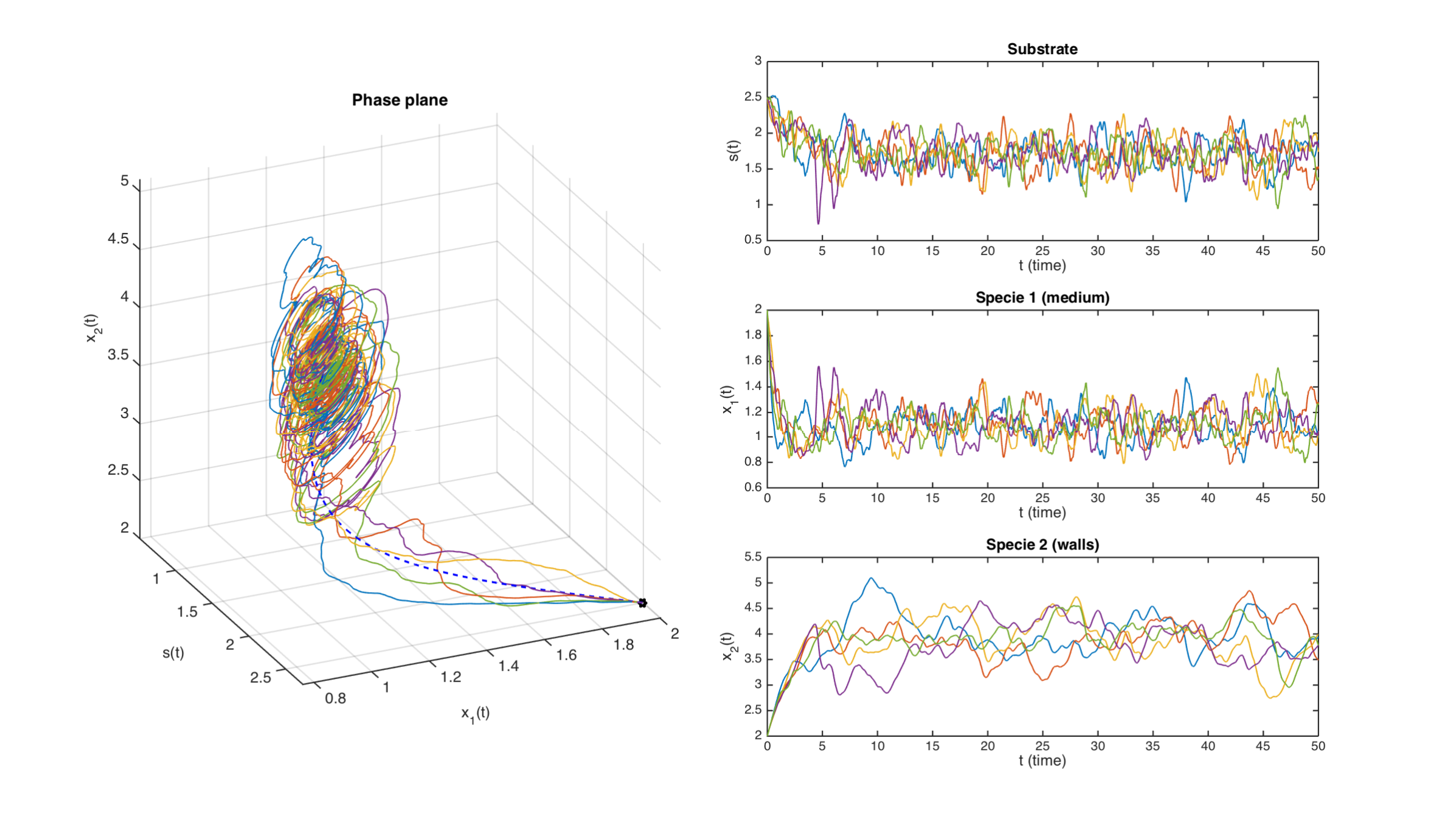}
\end{center}
\caption{Persistence of the species in the random chemostat model with wall growth (ii)}
\label{2_12}
\end{figure}

In the sequel, we only refer to the parameters to be changed respect to the last ones used and we suppose the rest to be the same than before. For instance, in Figure \ref{2_12} we will increase the quantity of noise to $\alpha=2$, the mean reversion constant to $\beta=4$ and the volatility constant to $\nu=0.7$, respect to the parameters used in the last figure.

\begin{figure}[H]
\begin{center}
\includegraphics[scale=0.27]{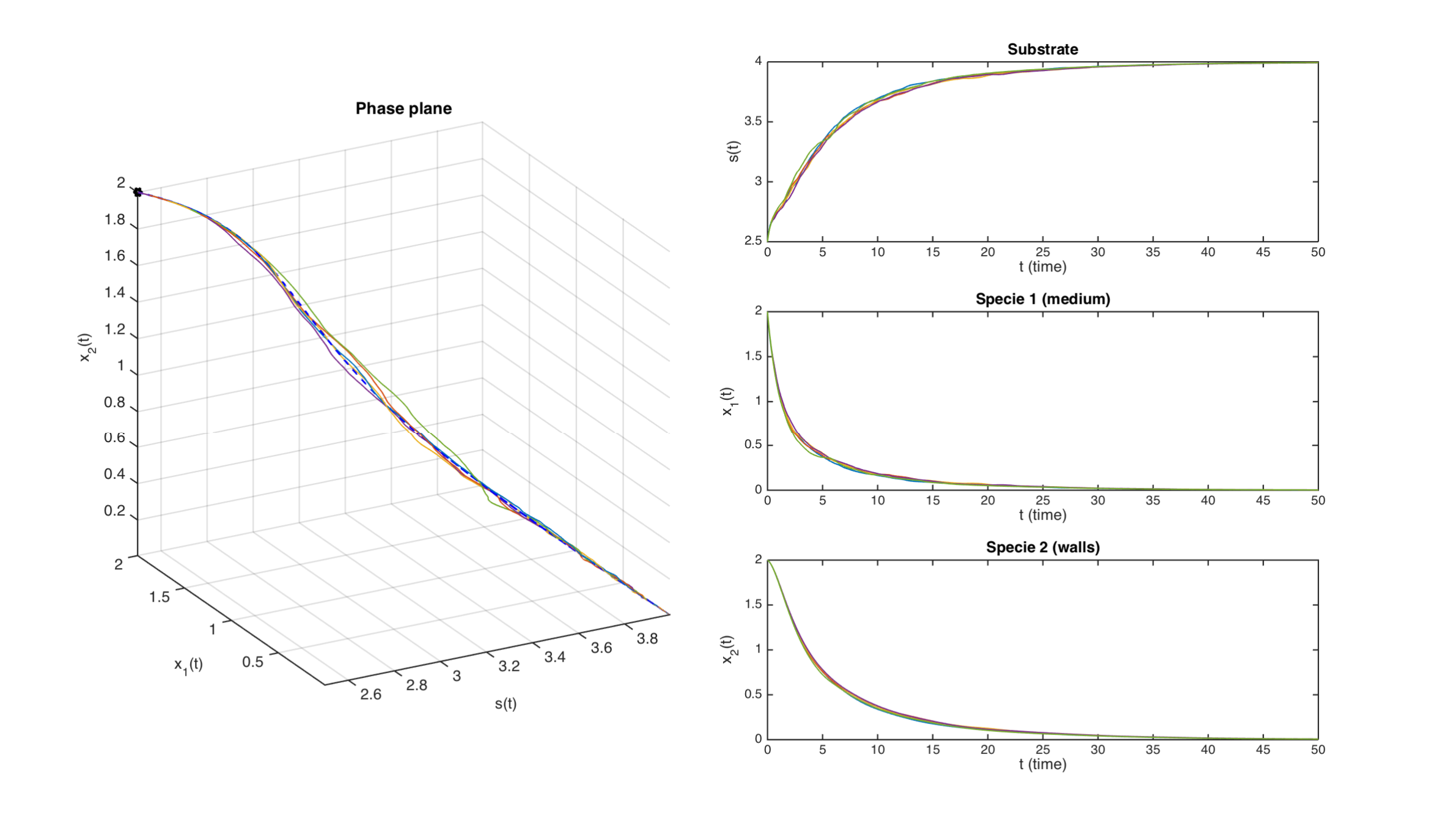}
\end{center}
\caption{Extinction of the species in the random chemostat model with wall growth}
\label{2_13}
\end{figure}

On another hand, in Figure \ref{2_13} we take $s_{in}=4$, $D=1.5$, $a=1.6$, $m=2$, $b=1$, $\nu=1.7$, $c=2.4$, $r_1=0.6$, $r_2=0.4$, $\alpha=0.5$, $\beta=1$, $\nu=0.2$ and we will take $s_0=2.5$, $x_{10}=2$, $x_{20}=2$ as initial values for the substrate and both species, respectively. Then, both species become extinct, as can be easily observed, which is quite logical in view of the values of the parameters, specially the collective death rate and the consumption rate of the species, which have been increased respect to the case of persistence.

\begin{figure}[H]
\begin{center}
\includegraphics[scale=0.27]{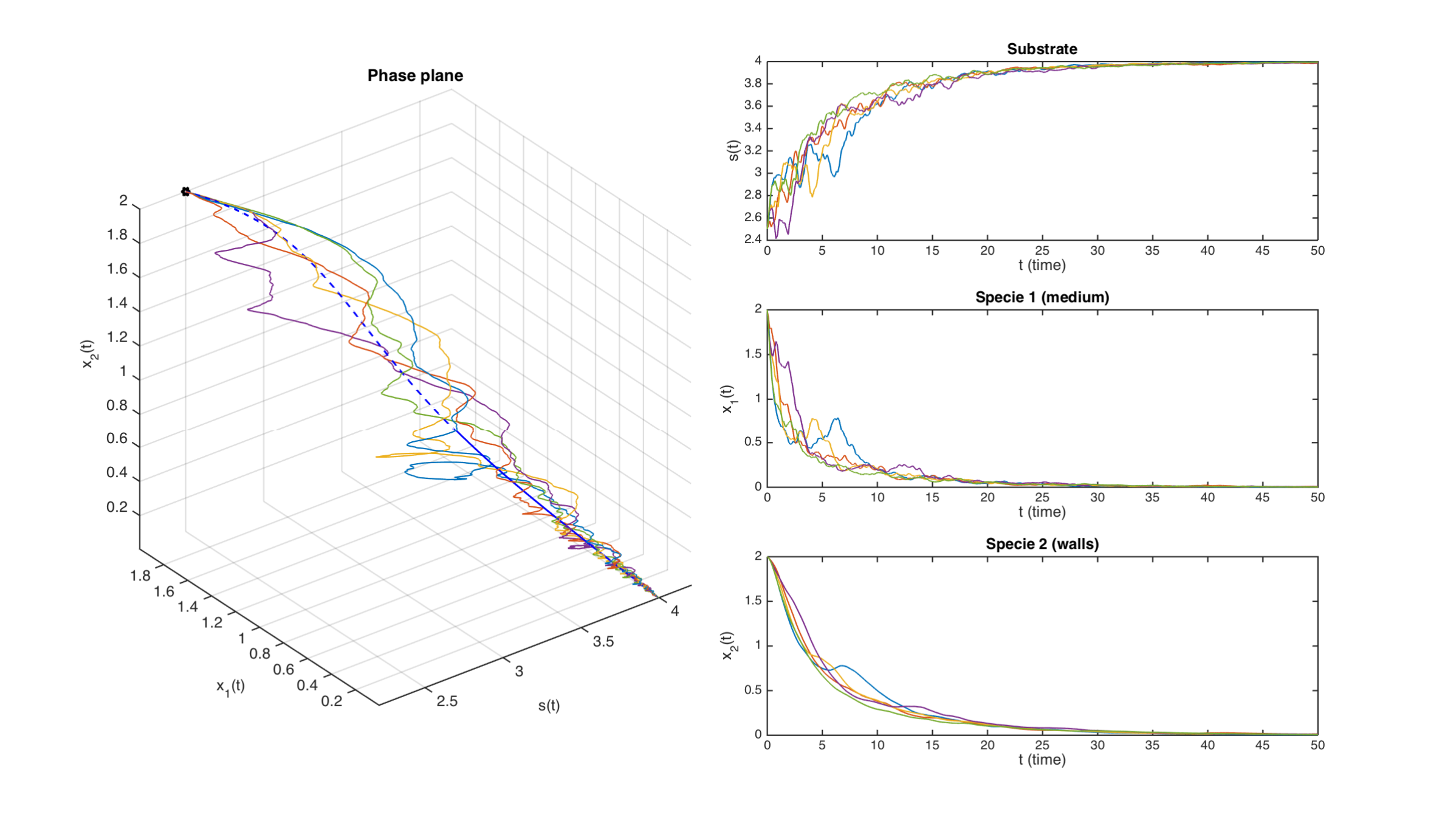}
\end{center}
\caption{Extinction of the species in the random chemostat model with wall growth (ii)}
\label{2_14}
\end{figure}

Finally, in Figure \ref{2_14} we increase the quantity of noise to $\alpha=2$, the mean reversion constant to $\beta=4$ and the volatility constant to $\nu=0.7$ and we remark that the rest of the parameters do not change respect to the last ones in Figure \ref{2_13}. In this case we can also see easily that both species become extinct, which is not surprising by taking into account what happened in the last figure and the new values of the parameters involved in the disturbances. 

\section{Stochastic chemostat model with wall growth}\label{s4}

In this section we analyze the chemostat model with wall growth \eqref{1}-\eqref{3} where the input flow is influenced by a standard Wiener process. Our main goal now is to study the effects produced by a non-bounded noise on a chemostat model with wall growth.\n

Let us recall the deterministic chemostat model with wall growth, which is given by the following differential system
\begin{eqnarray}
\f{ds}{dt}&=&D(s_{in}-s)-\f{ms}{a+s}x_1-\f{ms}{a+s}x_2+b\nu x_1,\label{dwall1}
\\[1.3ex]
\f{dx_1}{dt}&=&-(\nu+D)x_1+\f{cs}{a+s}x_1-r_1x_1+r_2x_2,\label{dwall2}
\\[1.3ex]
\f{dx_2}{dt}&=&-\nu x_2+\f{cs}{a+s}x_2+r_1x_1-r_2x_2,\label{dwall3}
\end{eqnarray}
\noindent where $s$, $x_1$ and $x_2$ denote concentrations of the nutrient and the two different microorganisms, respectively, and the parameters were introduced in Section \ref{intro}.\n

Our aim now is to perturb the input flow by the standard Wiener process in the same way that in the previous section, i.e., we replace the parameter $D$ by $D+\alpha\dot{\omega}(t)$ in \eqref{dwall1}-\eqref{dwall3}, where $\alpha>0$ represents the intensity of the noise and $\omega$ denotes the white noise, such that the following stochastic system understood in It\=o sense is obtained
\begin{eqnarray}
\nonumber
ds&=&\left[D(s_{in}-s)-\f{ms}{a+s}x_1-\f{ms}{a+s}x_2+b\nu x_1\right]dt+\alpha(s_{in}-s)d\omega(t),
\\[1.3ex]
\nonumber
dx_1&=&\left[-(\nu+D)x_1+\f{cs}{a+s}x_1-r_1x_1+r_2x_2\right]dt-\alpha x_1d\omega(t),
\\[1.3ex]
\nonumber
dx_2&=&\left[-\nu x_2+\f{cs}{a+s}x_2+r_1x_1-r_2x_2\right]dt.
\end{eqnarray}

After making use of the well-known conversion between It\=o and Stratonovich sense, we have the following stochastic system in Stratonovich sense
\begin{eqnarray}
\quad ds&=&\left[\bar{D}(s_{in}-s)-\f{ms}{a+s}x_1-\f{ms}{a+s}x_2+b\nu x_1\right]dt+\alpha(s_{in}-s)\circ d\omega(t),\label{swall1}
\\[1.3ex]
\quad dx_1&=&\left[-(\nu+\bar{D})x_1+\f{cs}{a+s}x_1-r_1x_1+r_2x_2\right]dt-\alpha x_1\circ d\omega(t),\label{swall2}
\\[1.3ex]
\quad dx_2&=&\left[-\nu x_2+\f{cs}{a+s}x_2+r_1x_1-r_2x_2\right]dt,\label{swall3}
\end{eqnarray}
\noindent where $\bar{D}=D+\f{\alpha^2}{2}$.

\subsection{Stochastic chemostat becomes a random chemostat}\label{s121}

In this section, we analyze the stochastic chemostat model with wall growth \eqref{swall1}-\eqref{swall3} by performing the variable change
\begin{eqnarray}
\sigma(t)&=&(s(t)-s_{in})e^{\alpha z^*(\theta_t\omega)},\label{sigma}
\\[1.3ex]
\kappa_1(t)&=&x_1(t)e^{\alpha z^*(\theta_t\omega)},\label{kappa1}
\\[1.3ex]
\kappa_2(t)&=&x_2(t),\label{kappa2}
\end{eqnarray}
\noindent where $z_*(\theta_t\omega)$ denotes the O-U process $z^*_{\beta,1}(\theta_t\omega)$ defined in Section \ref{s2}. The last variable $\kappa_2$ remains as $x_2$ due to the fact that equation \eqref{swall3} is not affected by the stochastic perturbation. For the sake of simplicity, we will write again $z^*$ instead of $z^*(\theta_t\omega)$ and $\sigma$, $\kappa_1$, $\kappa_2$ in place of $\sigma(t)$, $\kappa_1(t)$, $\kappa_2(t)$.\n

From \eqref{sigma}-\eqref{kappa2}, by differentiation, we obtain the following differential system satisfied by $\sigma$, $\kappa_1$ and $\kappa_2$, respectively,
{\small
\begin{eqnarray}
\quad \f{d\sigma}{dt}&=&-(\bar{D}+\alpha z^*)\sigma-\f{m(s_{in}+\sigma e^{-\alpha z^*})}{a+s_{in}+\sigma\ezm}\kappa_1-\f{m(s_{in}+\sigma\ezm)}{a+s_{in}+\sigma\ezm}\kappa_2\ez+b\nu\kappa_1,\label{rwall1}
\\[1.3ex]
\quad \f{d\kappa_1}{dt}&=&-(\nu+\bar{D}+\alpha z^*)\kappa_1+\f{c(s_{in}+\sigma\ezm)}{a+s_{in}+\sigma\ezm}\kappa_1-r_1\kappa_1+r_2\kappa_2\ez,\label{rwall2}
\\[1.3ex]
\quad \f{d\kappa_2}{dt}&=&-\nu\kappa_2+\f{c(s_{in}+\sigma\ezm)}{a+s_{in}+\sigma\ezm}\kappa_2+r_1\kappa_1\ezm-r_2\kappa_2.\label{rwall3}
\end{eqnarray}}

Now, we perform another variable change to transform the random system \eqref{rwall1}-\eqref{rwall3} into another one where the total biomass and the proportion of one of the species play an important and helpful role. To this end, we define the new state variables
\begin{eqnarray}
\kappa(t)&=&\kappa_1(t)+\kappa_2(t),\label{kappa}
\\[1.3ex]
\xi(t)&=&\f{\kappa_1(t)}{\kappa_1(t)+\kappa_2(t)}.\label{xi}
\end{eqnarray}

We also write in this case $\kappa$ and $\xi$ instead of $\kappa(t)$ and $\xi(t)$ in order to make the readability easier.\n

From \eqref{kappa}-\eqref{xi}, by differentiation, we have the following random system satisfied by $\sigma$, $\kappa$ and $\xi$, respectively,
{\small
\begin{eqnarray}
\nonumber
\f{d\sigma}{dt}&=&-(\bar{D}+\alpha z^*)\sigma-\f{m(s_{in}+\sigma\ezm)}{a+s_{in}+\sigma\ezm}\xi\kappa-\f{m(s_{in}+\sigma\ezm)}{a+s_{in}+\sigma\ezm}\ez\kappa(1-\xi)+b\nu\xi\kappa,\label{rpwall1}
\\[1.3ex]
\nonumber
\f{d\kappa}{dt}&=&-\nu\kappa-(\bar{D}+\alpha z^*)\kappa\xi+\f{c(s_{in}+\sigma\ezm)}{a+s_{in}+\sigma\ezm}\kappa+r_1\kappa\xi(\ezm-1)+r_2(1-\xi)\kappa(\ez-1),\label{rpwall2}
\\[1.3ex]
\nonumber
\f{d\xi}{dt}&=&-(\bar{D}+\alpha z^*)\xi(1-\xi)-r_1\xi+r_2\ez(1-\xi) -r_1(\ezm-1)\xi^2-r_2(\ez-1)\xi(1-\xi).\label{rpwall3}
\end{eqnarray}
}

Instead of analyzing now the existence and uniqueness of global solution of previous random system, we assume that there exists a unique global solution of the random chemostat model with wall growth. Particularly, from the random differential equations describing the dynamics of the substrate, we have the following equalities
\begin{eqnarray}
\nonumber
\left.\f{d\sigma}{dt}\right|_{\sigma=0}&=&-\f{ms_{in}}{a+s_{in}}\xi\kappa-\f{ms_{in}}{a+s_{in}}\ez\kappa(1-\xi)+b\nu\xi\kappa,
\\[1.3ex]
\nonumber
&=&-\f{ms_{in}}{a+s_{in}}\kappa\left[\xi+\ez(1-\xi)\right]
\\[1.3ex]
\nonumber
&=&\kappa\left[-\f{ms_{in}}{a+s_{in}}\left(\xi+\ez(1-\xi)\right)+b\nu\xi\right].\label{sigma0_wall}
\end{eqnarray}

Thus, $\sigma$ will remain positive as long as
\begin{equation*}
\f{ms_{in}}{a+s_{in}}\left(\xi(t)+\ezo(1-\xi(t))\right)\leq b\nu\xi(t)
\end{equation*}
\noindent holds true for every $t\geq 0$ and any $\omega\in\Omega$ or, equivalently, 
\begin{equation}
\ezo\leq\left(\f{b\nu\xi(t)(a+s_{in})}{ms_{in}}-\xi(t)\right)\f{1}{1-\xi(t)}.\label{cond_posit_wall_wiener}
\end{equation}

Thanks to \eqref{xi}, $0\leq \xi(t)\leq 1$ for every $t\geq 0$ since it is defined as the proportion of microorganisms in the liquid media. Hence, we have to distinguish the following situations:

\begin{itemize}
\item[$\blacksquare$] {\bf Case 1.- $\xi$ tends to one.} In this case, the quotient $1/(1-\xi(t))$ tends to infinity, thus \eqref{cond_posit_wall_wiener} could be true. Nevertheless, it would mean that $\kappa_2$ (or its corresponding state variable $x_2$) tends to zero. Consequently, the microbial biomass sticked on the walls of the culture vessel would become extinct, then we would recover the stochastic chemostat model without wall growth analyzed by Caraballo {\it et al} in \cite{CGL,corrigendumchapter}.\n
\item[$\blacksquare$] {\bf Case 2.- $\xi$ does not tend to one.} In this case, the right-hand side in \eqref{cond_posit_wall_wiener} is bounded for every $t\geq 0$ and any $\omega\in\Omega$, whereas $z^*$ could take arbitrary large values which means that \eqref{cond_posit_wall_wiener} does not hold true and then $\sigma$ can take negative values. In addition, it is not difficult to prove the following lower bound for $\sigma$, which is given by
\begin{equation}
\sigma(t)>-(a+s_{in})\ezo,\label{sn}
\end{equation}
\noindent similarly to the proof made in \cite{corrigendumchapter} (Theorem ), that implies
\begin{equation}
s(t)>-a.\label{nn}
\end{equation}

As explained in \cite{CGL,corrigendumchapter}, it is not realistic at all from the biological point of view and this is the main reason because we will not develop a deeper analysis of this kind of noise on the chemostat model with wall growth.
\end{itemize}

\subsection{Numerical simulations and final comments}\label{s122}

In this section we show some numerical simulations to remark the drawbacks found when using the Wiener process to perturb the input flow in the chemostat model with wall growth.\n

In each figure some different panels are displayed. On the one hand, there is a big panel on the left, where we can see the phase plane showing the dynamics of the stochastic chemostat model with wall growth in Stratonovich sense \eqref{swall1}-\eqref{swall3}. On the other hand, three different panels are stated on the right to describe the dynamics of the substrate and both species individually, depending on the time.\n

The blue dashed lines in the big panel represent the solutions of the deterministic system \eqref{dwall1}-\eqref{dwall3}. In addition, we set $s_{in}=1$, $a=0.6$, $m=3$, $r_1=0.6$, $r_2=0.4$ and we consider $(s(0),x_1(0),x_2(0))=(5,2.5,2.5)$ as initial data. In this way, we present different cases where the value of the rest of the parameters change and the intensity of the noise increases or decreases, in order to show the effect of each one on the dynamics of our model.\newpage

\begin{figure}[H]
\begin{center}
\includegraphics[scale=0.35]{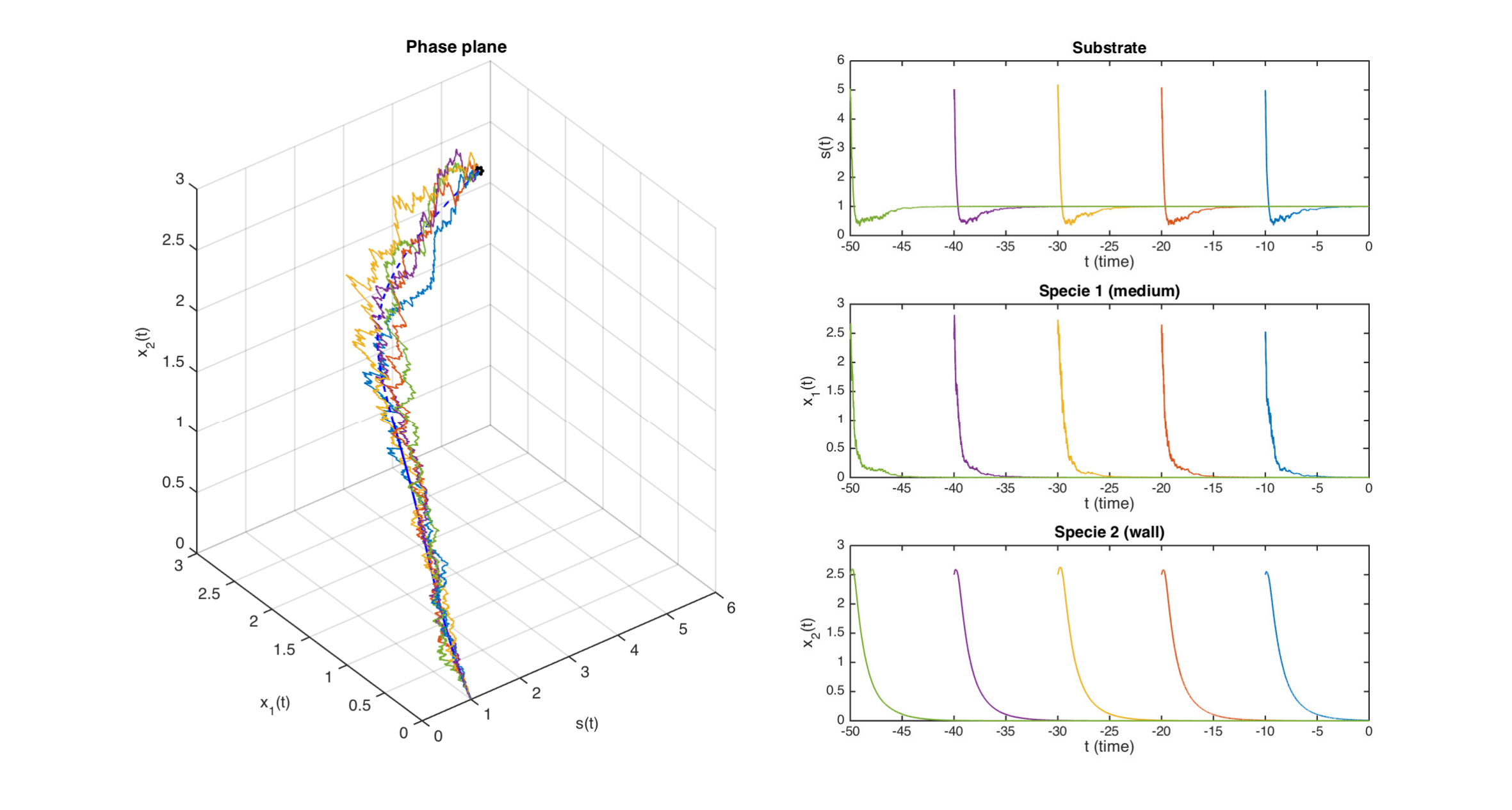}
\caption{Extinction of both species. $\alpha=0.5$}
\label{sim6}
\end{center}
\end{figure}

On the one hand, in Figure \ref{sim6} we take $D=3$, $b=2$, $\nu=0.2$, $c=1.5$ and we choose $\alpha=0.5$. We can see that every state variable remains in the first octant in this case, i.e., they are all positive, however both species, the one in the medium or liquid media and also the one sticked on the walls, become extinct.

\begin{figure}[H]
\begin{center}
\includegraphics[scale=0.35]{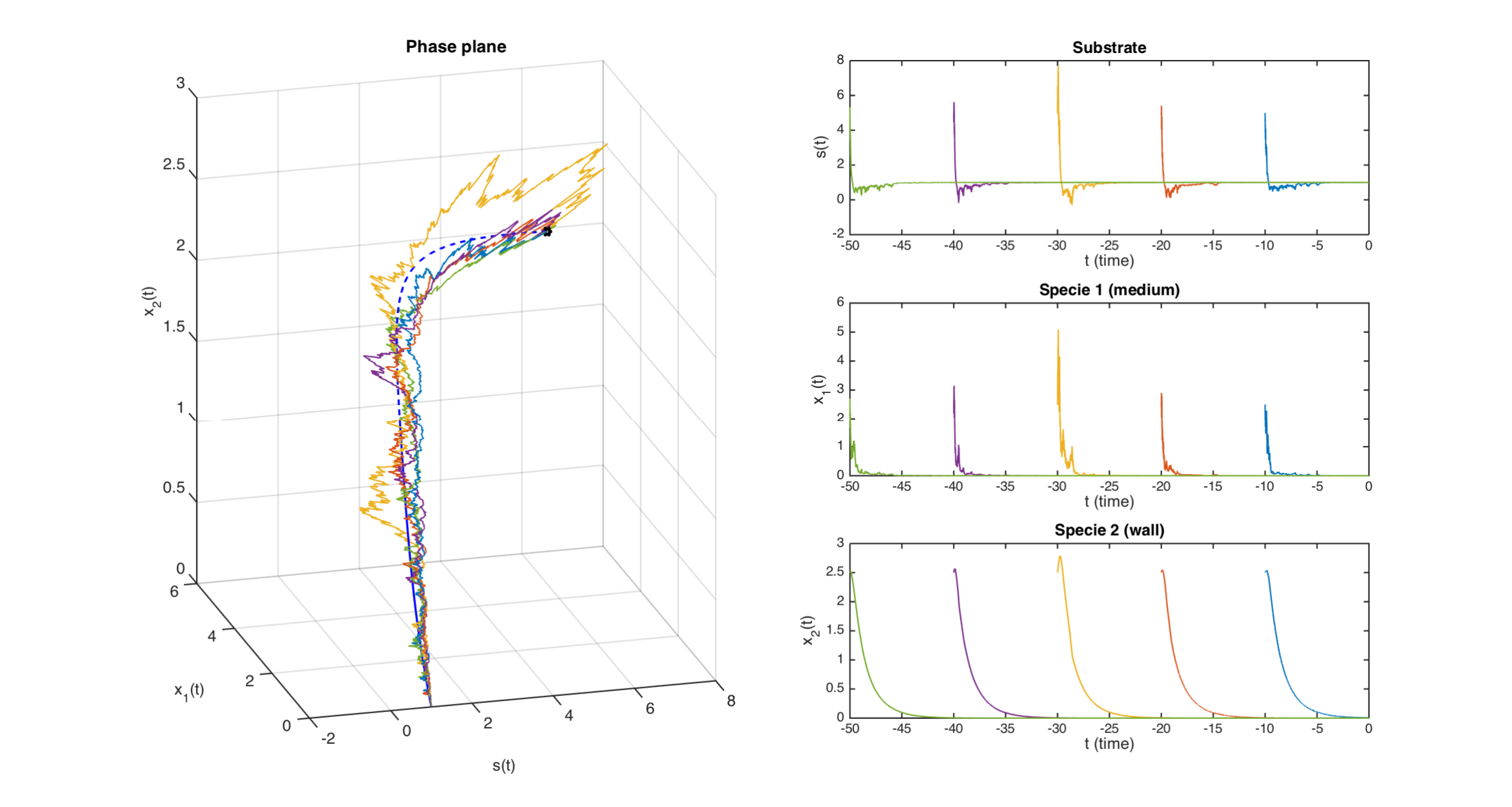}
\caption{Extinction of both species. $\alpha=1.5$}
\label{sim7}
\end{center}
\end{figure}

Similarly, in Figure \ref{sim7} we take $D=3$, $b=2$, $\nu=0.2$, $c=1.5$ and we increase the intensity of the noise to $\alpha=1.5$. In this case, both species become extinct and, moreover, the substrate reaches negative values for some times which is totally unrealistic from the biological point of view, as explained previously.

\begin{figure}[H]
\begin{center}
\includegraphics[scale=0.35]{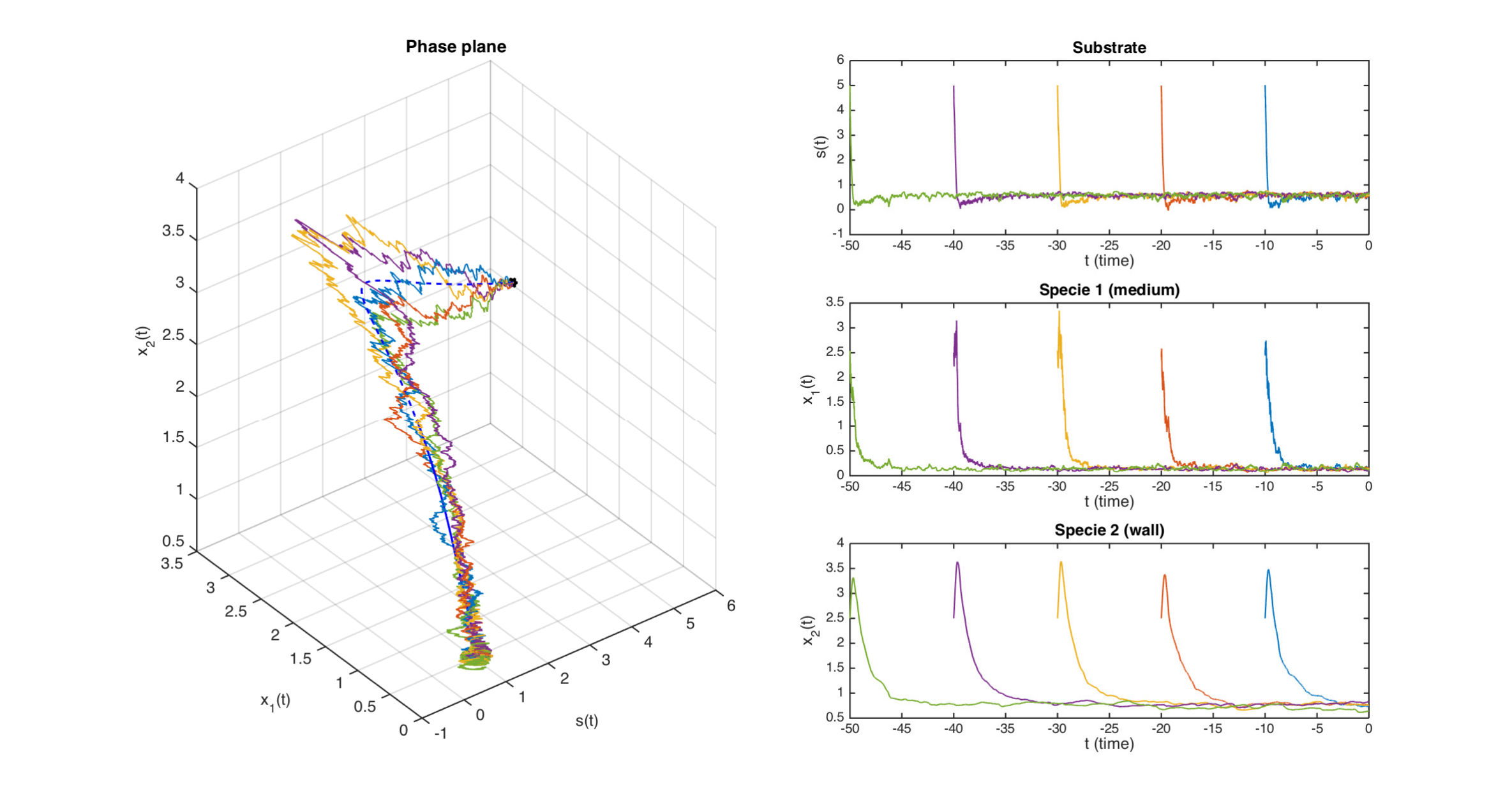}
\caption{Persistence of both species. $\alpha=0.5$}
\label{sim8}
\end{center}
\end{figure}

Now, in Figure \ref{sim8} we take $D=3$, $b=0.5$, $\nu=1.2$ and $\alpha=0.5$. In this case, it is not difficult to notice that both species persist. Nevertheless, the dynamics of the substrate clearly cross the line $s=0$ which is an important inconvenient.

\newpage

\begin{figure}[H]
\begin{center}
\includegraphics[scale=0.35]{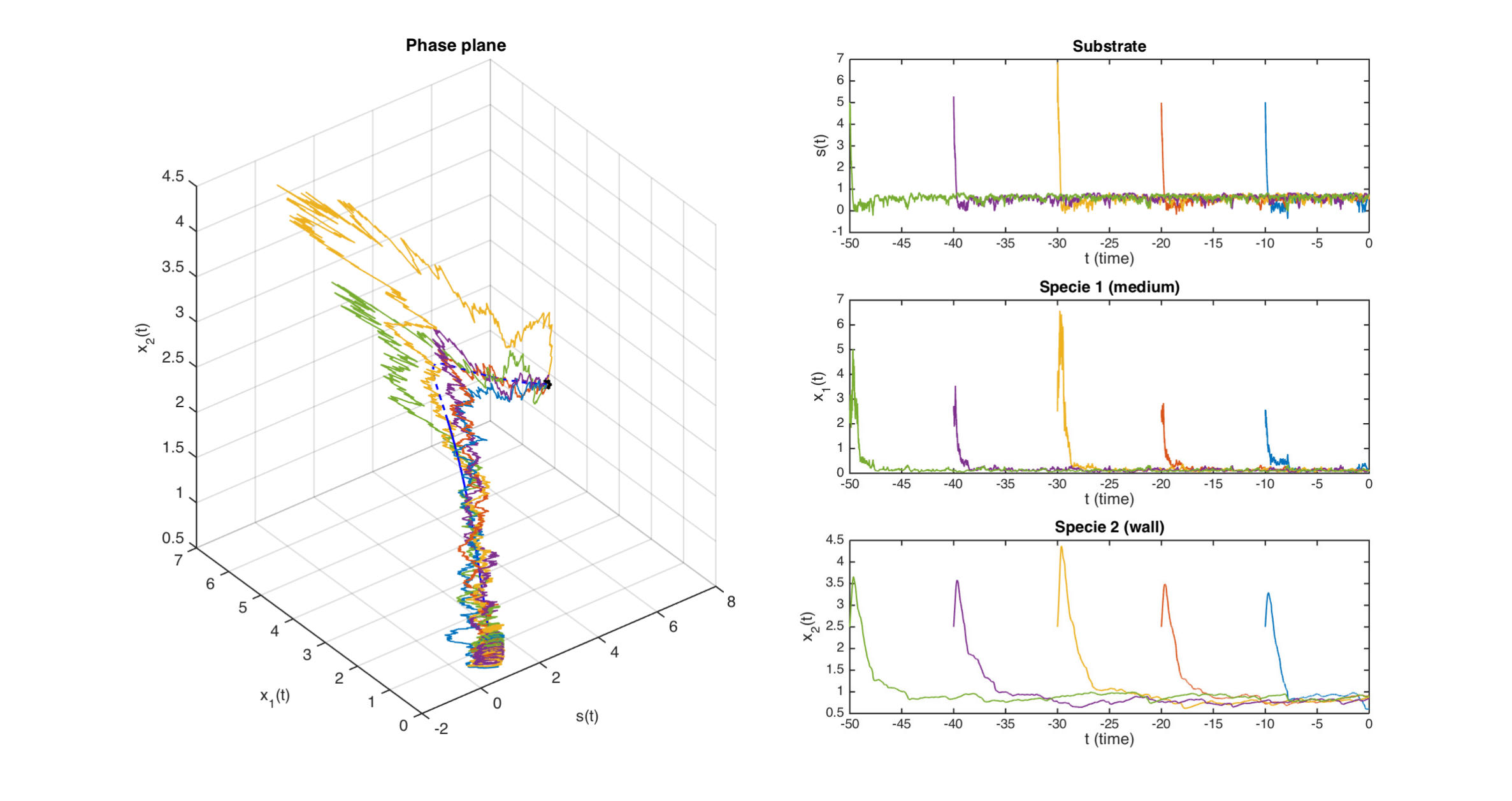}
\caption{Persistence and extinction. $\alpha=1.5$}
\label{sim9}
\end{center}
\end{figure}

Finally, in Figure \ref{sim9} we take $D=3$, $b=0.5$, $\nu=1.2$ and $\alpha=1.5$. In this case, the species which are sticked on the walls of the culture vessel persist whereas the ones in the liquid media become extinct. In addition, we also find some drawbacks since the substrate reaches negative values.

\section{Comparison between random and stochastic models} \label{s5}

In this section we make a comparison between both ways of modeling randomness and stochasticity on the input flow in the chemostat model with wall growth \eqref{rp1}-\eqref{rp3}  in order to show the differences between the effects caused by the usual standard Wiener process and remark the important advantages of using a bounded noise, in particular, the Ornstein-Uhlenbeck process introduced in Section \ref{s2}.\n

From the theoretical point of view, we highlight the following points

\begin{itemize}
\item {\bf Stochastic process:} the main difference between Sections \ref{s3} and \ref{s4} is clearly the stochastic process used to model the fluctuations on the input flow in the chemostat model. Firstly, in Section \ref{s3} we used the suitable O-U process introduced in Section \ref{s2} to perturb the input flow by that bounded noise obtaining the random chemostat model (RCM) with wall growth \eqref{1}-\eqref{3}. Therefore, in Section \ref{s4}we use the common standard Wiener process to perturb the input flow and we obtain the stochastic chemostat model (SCM) \eqref{1}-\eqref{3}.\n
\item {\bf Positiveness of solution:} we would like to remark that in the (RCM) it is possible to prove that every solution remains positive for any positive initial value (see Theorem \ref{euap}) whereas in the (SCM) it is not possible to prove it due to the huge fluctuations (positive or negative) that the white noise could present (see \eqref{sn} and \eqref{nn}).\n
\item {\bf Persistence and coexistence of the species:} we notice firstly that this point is the most important and interesting for biologists. In the (RCM) the persistence of the total microbial biomass, $x=x_1+x_2$, is proved as well as the coexistence of both species individually, by assuming the following conditions
\begin{equation}
\nu<c-D\xi^*_l\quad\quad\text{and}\quad\quad \nu+b_2<\f{z^*_l}{a+\f{z^*_u}{c}},\label{cp}
\end{equation}
\noindent where $\xi^*_l$, $z^*_l$ and $z^*_u$ are deterministic constants defined by \eqref{boundsxiat} and \eqref{ze}, respectively.\n

We would like to highlight that both conditions in \eqref{cp} essentially represent some restrictions on the dilution rate, on the disturbances on the input flow and also on the death collective rate, which is totally logical from the biological point of view.\n

In particular, if the dilution rate, or its equivalent input flow, were too large, then the microbial biomass would not be able to have access to the nutrient which would mean the extinction of both species and, furthermore, much more quantity of microbial biomass would be removed from the culture vessel to the collection vessel which would also increase significantly the probability of the extinction. In addition, the disturbances on the input flow cannot be too large since we want to avoid the drawbacks found when modeling the disturbances by means of the white noise, in fact, it is not possible to provide conditions to ensure the persistence or coexistence of the species in the (SCM).\n

Apart from that, it is not surprising the presence of the death collective rate in both conditions \eqref{cp} since it must be difficult to prove the persistence of both species if this parameter is too large. Therefore, the conditions required to get the persistence of both species are, as already pointed out, absolutely reasonable from the biological point of view.\n
\item {\bf Theoretical framework:} this is also an important point between the (RCM) and the (SCM). In the last one, a framework based on the theory of random dynamical systems and pullback attractors is needed (see \cite{CGL,corrigendumchapter}) although we did not present this framework in this work due to the fact that we just wanted to show the drawbacks found when using the Wiener process and we tried to avoid the analysis since the resulting (SCM) was not realistic from the biological point of view.\n

In addition, every result based on the theory of random dynamical systems is proved by using the pullback convergence which is not as natural as the forwards one used when analyzing the (RCM). Moreover, the pullback convergence do not provide enough information about the long-time behavior of the models in many situations, which is another inconvenient of using the white noise.
\end{itemize}

On the other hand, we would like to present some numerical simulations where a typical realization  of the (SCM) analyzed in Section \ref{s4} and two typical ones of the (RCM) analyzed in Section \ref{s3} will be plotted together in order to see easily the differences of modeling the disturbances on the input flow in the chemostat model by using both the white noise (orange lines) and the Ornstein-Uhlenbeck process (red and green lines).\n

\begin{figure}[H]
\begin{center}
\includegraphics[scale=0.35]{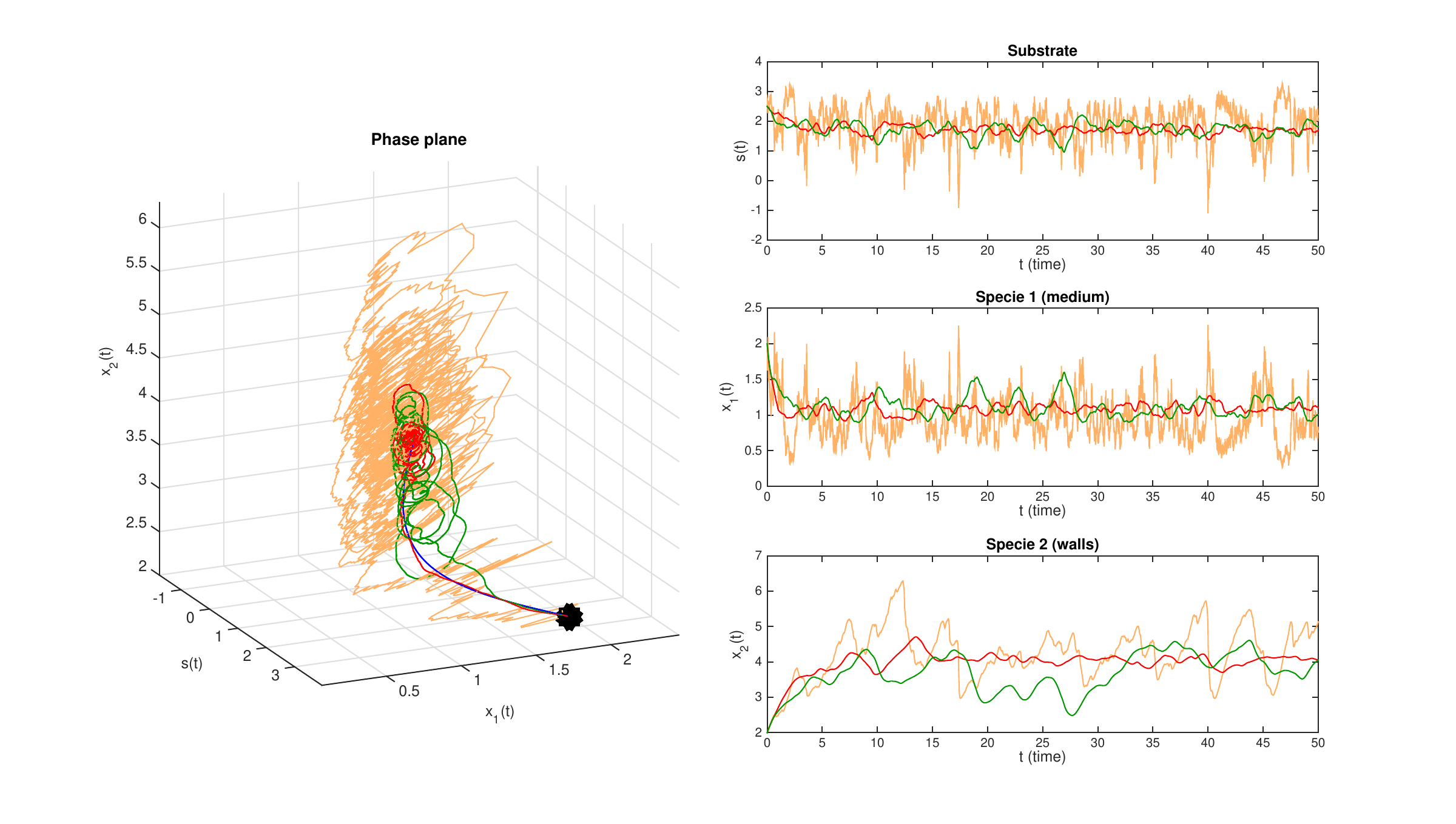}
\end{center}
\caption{Comparison in case of persistence}
\label{2_15}
\end{figure}

Firstly, in Figure \ref{2_15} we take $s_{in}=4$, $D=2$, $a=1.6$, $m=2$, $b=0.5$, $\nu=1.2$, $c=3$, $r_1=0.2$, $r_2=0.4$ and we consider $s_0=2.5$, $x_{10}=2$, $x_{20}=2$ as initial values for the substrate and both species, respectively.  We also choose $\alpha=0.8$, $\nu=0.7$ and both $\beta=2$ (red line) and $\beta=1$ (green line). In this case, we can see that both species persist and we remark the huge disturbances obtained in case of using the white noise respect to the Ornstein-Uhlenbeck process, in fact, these disturbances can be also observed to affect to the specie $x_2$ even though it is not affected by the random input flow directly (see \eqref{rp3}), which do not happen when perturbing the input flow by means of the bounded noise.

\begin{figure}[H]
\begin{center}
\includegraphics[scale=0.35]{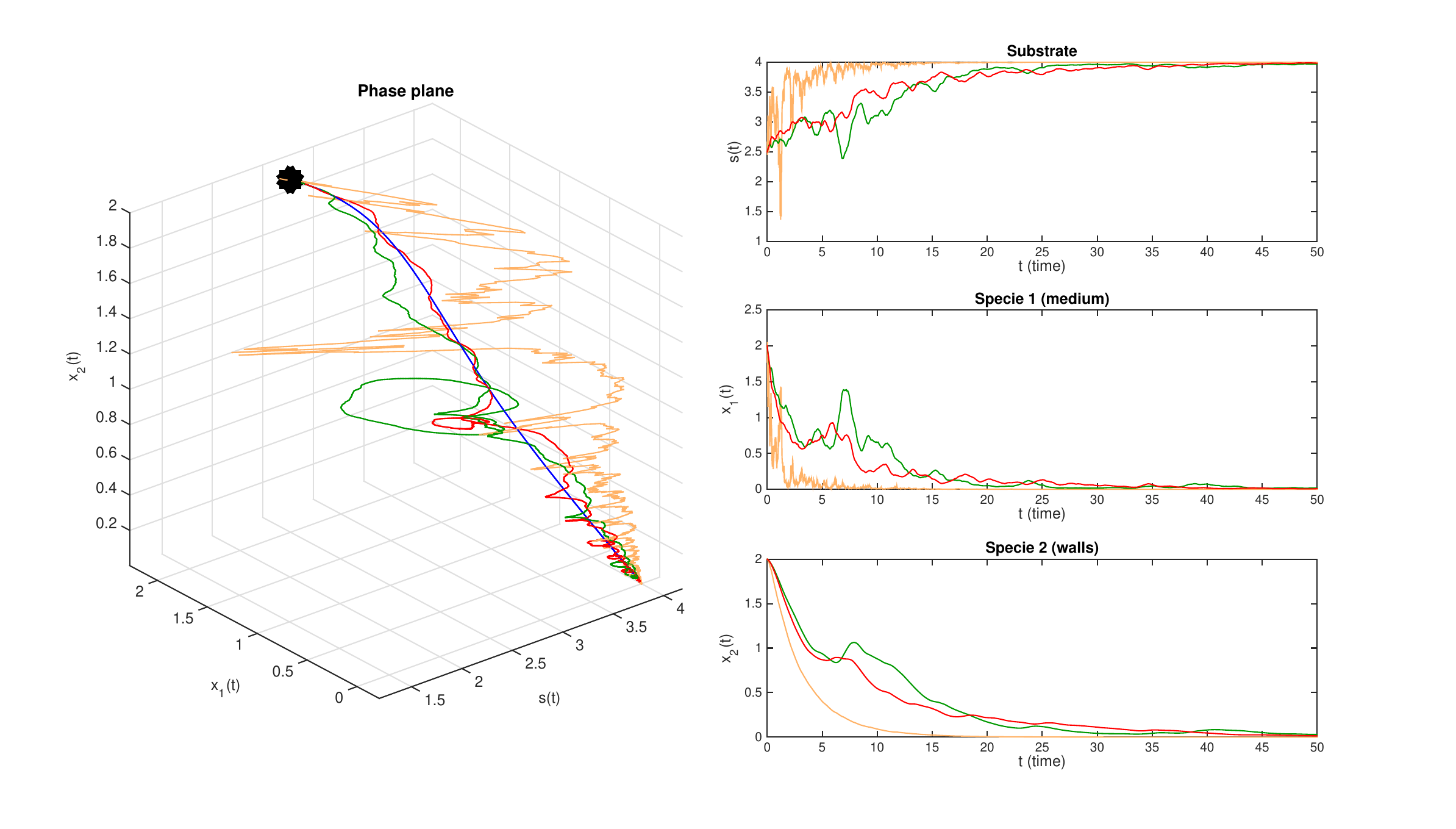}
\end{center}
\caption{Comparison in case of extinction}
\label{2_16}
\end{figure}

Finally, in Figure \ref{2_16} we take $s_{in}=4$, $D=1.5$, $a=1.6$, $m=2$, $b=1$, $\nu=1.7$, $c=2.4$, $r_1=0.6$, $r_2=0.4$ and we consider $s_0=2.5$, $x_{10}=2$, $x_{20}=2$ as initial values for the substrate and both species, respectively. In this case we increase the quantity of noise, respect to the last case, to $\alpha=1.5$, $\nu=0.7$ and both $\beta=2$ (red) and $\beta=1$ (green). We can observe that both species become extinct and we remark again the significant disturbances when using the white noise respect to the case when considering the Ornstein-Uhlenbeck process.

\section{Final comments}\label{s6}

In this paper we prove again the important advantages of using the suitable O-U process introduced in Section \ref{s2} when modeling fluctuations in chemostat models which, in addition, allow us to avoid the important drawbacks found when using the usual white noise.\n

In this case we can go further than in the simplest chemostat model in \cite{CGLR} since now we have two different species and we are able to prove not only the persistence of the total biomass but also the coexistence of both species which is much more interesting from the biological point of view and gives our work much more value.\n

Apart from that, due to the fact that most of fluctuations in the real life are bounded, the O-U process used in this work can be also used to model randomness in other systems in life science such as the logistic equation, Lotka-Volterra systems, SIR model... where interesting results are also expected to be achieved.


\bibliographystyle{aims.bst}

\bibliography{submitted}

\providecommand{\href}[2]{#2}
\providecommand{\arxiv}[1]{\href{http://arxiv.org/abs/#1}{arXiv:#1}}
\providecommand{\url}[1]{\texttt{#1}}
\providecommand{\urlprefix}{URL }
\begin{thebibliography}{10}

\bibitem{alazzawi2017}
\newblock S.~Al-\uppercase{a}zzawi, J.~Liu and X.~Liu,
\newblock Convergence rate of synchronization of systems with additive noise,
\newblock \emph{Discrete Contin. Dyn. Syst. Ser. B}, \textbf{22} (2017),
  227--245.

\bibitem{arnold}
\newblock L.~Arnold,
\newblock \emph{Random Dynamical Systems},
\newblock Springer Berlin Heidelberg, 1998.

\bibitem{CGL}
\newblock T.~Caraballo, M.~J. {Garrido-\uppercase{a}tienza} and
  J.~{{L{\'{o}}pez-de-la-\uppercase{c}ruz}},
\newblock \emph{Some Aspects Concerning the Dynamics of Stochastic Chemostats},
  vol.~69, chapter~11, 227--246,
\newblock Springer International Publishing, Cham, 2016.

\bibitem{CGLii}
\newblock T.~Caraballo, M.~J. {Garrido-\uppercase{a}tienza} and
  J.~{{L{\'{o}}pez-de-la-\uppercase{c}ruz}},
\newblock Dynamics of some stochastic chemostat models with multiplicative
  noise,
\newblock \emph{Communications on Pure and Applied Analysis}, \textbf{16}
  (2017), 1893--1914.

\bibitem{corrigendumchapter}
\newblock T.~Caraballo, M.~J. {Garrido-\uppercase{a}tienza},
  J.~{{L{\'{o}}pez-de-la-\uppercase{c}ruz}} and A.~Rapaport,
\newblock Corrigendum to "\uppercase{S}ome aspects concerning the dynamics of
  stochastic chemostats", 2017,
\newblock \emph{arXiv:1710.00774 [math.DS]}.

\bibitem{CGLR}
\newblock T.~Caraballo, M.~J. {Garrido-\uppercase{a}tienza},
  J.~{{L{\'{o}}pez-de-la-\uppercase{c}ruz}} and A.~Rapaport,
\newblock Modeling and analysis of random and stochastic input flows in the
  chemostat model,
\newblock \emph{Discrete \& Continuous Dynamical Systems - Series B.}

\bibitem{caraballo-book}
\newblock T.~Caraballo and X.~Han,
\newblock \emph{Applied Nonautonomous and Random Dynamical Systems, Applied
  Dynamical Systems},
\newblock Springer International Publishing, 2016.

\bibitem{CHKR}
\newblock T.~Caraballo, X.~Han, P.~{E. Kloeden} and A.~Rapaport,
\newblock \emph{Continuous and Distributed Systems II}, chapter Dynamics of Non
  autonomous Chemostat Models, 103--120,
\newblock Springer International Publishing, Cham, 2015.

\bibitem{caraballo4}
\newblock T.~Caraballo, X.~Han and P.~E. Kloeden,
\newblock Chemostats with time-dependent inputs and wall growth,
\newblock \emph{Applied Mathematics and Information Sciences}, \textbf{9}
  (2015), 2283--2296.

\bibitem{caraballo1}
\newblock T.~Caraballo, P.~E. Kloeden and B.~Schmalfuss,
\newblock Exponentially stable stationary solutions for stochastic evolution
  equations and their perturbation,
\newblock \emph{Applied Mathematics and Optimization}, \textbf{50} (2004),
  183--207.

\bibitem{MR}
\newblock M.~{El Hajji} and A.~Rapaport,
\newblock Practical coexistence of two species in the chemostat - a slow-fast
  characterization,
\newblock \emph{Mathematical Biosciences}, \textbf{218} (2009), 33--39.

\bibitem{HLRS17}
\newblock J.~Harmand, C.~Lobry, A.~Rapaport and T.~Sari,
\newblock \emph{The Che\-mos\-tat: Mathematical Theory of Micro-organisms
  Cultures},
\newblock Wiley, Chemical Engineering Series, John Wiley {\&} Sons, Inc., 2017.

\bibitem{imhof-walcher}
\newblock L.~Imhof and S.~Walcher,
\newblock Exclusion and persistence in deterministic and stochastic chemostat
  models,
\newblock \emph{Journal of Differential Equations}, \textbf{217} (2005),
  26--53.

\bibitem{J74}
\newblock H.~W. Jannasch,
\newblock Steady state and the chemostat in ecology,
\newblock \emph{Limnology and Oceanography}, \textbf{19} (1974), 716--720.

\bibitem{HS}
\newblock V.~{Sree Hari Rao} and P.~{Raja Sekhara Rao},
\newblock \emph{Dynamic Models and Control of Biological Systems},
\newblock Springer-Verlag, Heidelberg, 2009.

\bibitem{XYZ13}
\newblock C.~Xu, S.~Yuan and T.~Zhang,
\newblock Asymptotic behavior of a che\-mos\-tat model with stochastic
  perturbation on the dilution rate,
\newblock \emph{Abstract and Applied Analysis}, 1--11.

\end{thebibliography}

\end{document}